\documentclass[english]{article}
\usepackage[T1]{fontenc}
\usepackage[latin9]{inputenc}
\usepackage{geometry}
\usepackage{amsmath}
\usepackage{amssymb}
\usepackage{graphicx}
\usepackage{esint}

\makeatletter

\providecommand{\tabularnewline}{\\}

\makeatother

\usepackage{babel}
\begin{document}

\title{A new multigroup method for cross-sections that vary rapidly in energy}

\author{T. S. Haut and C. Ahrens and A. Jonko and R. Lowrie and A. Till \\
Los Alamos National Laboratory\\
Los Alamos, NM 87545}
\maketitle
\begin{abstract}
We present a numerical method for solving the time-independent thermal
radiative transfer (TRT) equation or the neutron transport (NT) equation
when the opacity or cross-section varies rapidly in energy (frequency)
on the microscale $\varepsilon$; $\varepsilon$ corresponds to the
characteristic spacing between absorption lines or resonances, and
is much smaller than the macroscopic energy (frequency) variation
of interest. The approach is based on a rigorous homogenization of
the TRT/NT equation in the energy (frequency) variable. Discretization
of the homogenized TRT/NT equation results in a multigroup-type system,
and can therefore be solved by standard methods.

We demonstrate the accuracy and efficiency of the approach on three
model problems. First we consider the Elsasser band model with constant
temperature and a line spacing $\varepsilon=10^{-4}$. Second, we
consider a neutron transport application for fast neutrons incident
on iron, where the characteristic resonance spacing $\varepsilon$
necessitates $\approx16,000$ energy discretization parameters if
Planck-weighted cross sections are used. Third, we consider an atmospheric
TRT problem for an opacity corresponding to water vapor over a frequency
range $1000-2000\,\mathrm{cm}^{-1}$, where we take $12$ homogeneous
layers between $1\,$km - $15\,$km, and temperature/pressure values
in each layer from the standard US atmosphere. For all three problems,
we demonstrate that we can achieve between $0.1$ and $1$ percent
relative error in the solution, and with several orders of magnitude
fewer parameters than a standard multigroup formulation using Planck-weighted
opacities for a comparable accuracy.
\end{abstract}

\section{Background}

Thermal radiative transfer (TRT) plays a key role in a number of scientific
and engineering disciplines. For example, resolving the radiation
field in three-dimensional cloudy atmospheres is key to understanding
a number of atmospheric science and remote sensing problems \cite{MAR-DAV:2005}.
In many TRT problems, there is rapid variation in the opacity with
energy (or frequency) due to bound-bound and bound-free transitions.
In fact, for broad-band TRT problems there can be hundreds of thousands
of absorption lines, whose widths are many times smaller than the
overall energy range of interest. This fine scale structure in the
opacities, coupled with discretizing the spatial and angular variables,
places large demands on computational resources. Thus, researchers
have sought methods to ``average'' or ``homogenize'' the opacities
and derive so-called ``grey'' or frequency independent approximations,
thereby reducing the complexity of solving the full TRT problem. Many
opacity homogenization techniques have been developed over the years.
Here we mention only those most closely related with the method developed
in this paper. 

A commonly used averaging technique (the multigroup or picket fence)
starts by integrating the transport equation over the energy interval
$\left[E_{g},E_{g+1}\right]$. This formerly results in a transport
equation for the group averaged intensity,
\[
\psi_{g}\left(\mathbf{x},\boldsymbol{\Omega}\right)=\int_{E_{g}}^{E_{g+1}}\psi\left(\mathbf{x},\boldsymbol{\Omega},E\right)dE.
\]
However, the group averaged opacities $\sigma_{g}$ in the resulting
multigroup equations depend on the unknown solution $\psi\left(\mathbf{x},\boldsymbol{\Omega},E\right)$,
and an approximation is therefore needed to close the system. Typically,
either a Rosseland or Plank mean opacity is used; generally, these
closures are only accurate when $\left[E_{g},E_{g+1}\right]$ is small
relative to the variation of $\sigma\left(E\right)$ or in certain
limiting physical regimes (e.g. the optically thick limit). See, e.g.,
\cite{POMRAN:1971} for details. We remark that the averaging method
developed in this paper does not require one to postulate such a closure
relationship. 

The multiband method \cite{CUL-POM:1980} is another approach for
averaging the transport equation. The formulation results in a multigroup-type
system for
\begin{equation}
\psi_{g,b}\left(\mathbf{x},\boldsymbol{\Omega}\right)=\int_{\sigma_{g,b}}^{\sigma_{g,b+1}}\int_{E_{g}}^{E_{g+1}}\delta\left(\sigma\left(E\right)-\xi\right)\psi\left(\mathbf{x},\boldsymbol{\Omega},E\right)dEd\xi.\label{eq:multiband average}
\end{equation}
The multiband equations depend on an averaged opacity $\sigma_{g,b}$
that again involves the unknown solution $\psi\left(\mathbf{x},\boldsymbol{\Omega},E\right)$.
Like the multigroup method, this requires one to make some approximation
in order to close the system. We note that the group averaged solution
is recovered from (\ref{eq:multiband average}) via summing over bands
$b$ within each group $g$.

In the context of atmospheric TRT calculations, the correlated k-Distribution
method \cite{ARK-GRO:1972} and the multigroup k-Distribution method
\cite{MOD-ZHA:2001} have been shown to drastically reduce the computational
cost of direct line-by-line calculations, and have similarities to
the current homogenization approach. Another method that is related
to the current homogenization approach is the so-called Opacity Distribution
Function (ODF) method \cite{STR-KUR:1966}, which also takes a statistical
approach toward coarse-graining the TRT equation in energy. See \cite{MODEST:2013}
for a clear overview of these, and related, methods.

Finally, let us mention that the homogenization in energy of transport
equations in the absence of scattering can be analyzed using techniques
developed by Tartar \cite{TARTAR:1990}. Unlike when the opacity has
rapid spatial variation (see \cite{DUM-GOL:2000}), the homogenized
version of the transport equation \textit{is not }simply another transport
equation with a homogenized opacity. In fact, the homogenized equation
is an integro-differential equation, where the integral equation is
nonlocal in the spatial variable. An alternative homogenization approach
for solving the time-dependent TRT equation, in the absence of scattering,
has also been pursued in \cite{ET-GO-SA:2010} and \cite{MAT-SAL:2013},
where the authors develop homogenized equations on an enlarged phase
space. This approach is similar in spirit to that taken in the current
paper; however, one advantage of that given here is the ability to
handle scattering (in angle and energy), as well as to reduce the
numerical computation to a standard multigroup formulation. Finally,
we remark that, when the opacity is of the form $\sigma_{\varepsilon}\left(E,T\right)=\sigma_{0}\left(E/\varepsilon,E,T\right)$,
where $\sigma_{0}\left(\kappa,E,T\right)$ is almost periodic in $\kappa$,
then (\ref{eq:transport equation on extended phase space, intro})
and (\ref{eq:psi_0 as weighting against Young measure, intro}) is
analogous to the two-scale homogenization theory developed in \cite{ALLAIR:1992}.

\section{Outline of the new multigroup method}

Here we outline a computational method for the transport equation
\begin{eqnarray}
\boldsymbol{\Omega}\cdot\nabla_{\mathbf{x}}\psi_{\varepsilon}+\sigma_{\varepsilon}\left(T\left(\mathbf{x}\right),E\right)\psi_{\varepsilon} & = & \sigma_{\varepsilon}^{a}\left(T\left(\mathbf{x}\right),E\right)S\left(\mathbf{x},E\right)+Q_{0}\left(\mathbf{x},\boldsymbol{\Omega},E\right)+\nonumber \\
 &  & \int_{\mathbb{S}^{2}}\int_{0}^{\infty}\Sigma^{s}\left(\mathbf{x},E,E',\boldsymbol{\Omega}\cdot\boldsymbol{\Omega}'\right)\psi_{\varepsilon}dE'd\boldsymbol{\Omega}',\label{eq:transport equation, intro}
\end{eqnarray}
where 
\[
\Sigma^{s}\left(\mathbf{x},E,E',\boldsymbol{\Omega}\cdot\boldsymbol{\Omega}'\right)=\sigma^{s}\left(T\left(\mathbf{x}\right),E'\right)K\left(E,E',\boldsymbol{\Omega}\cdot\boldsymbol{\Omega}'\right),
\]
\[
\sigma_{\varepsilon}\left(T\left(\mathbf{x}\right),E\right)=\sigma_{\varepsilon}^{a}\left(T\left(\mathbf{x}\right),E\right)+\sigma^{s}\left(T\left(\mathbf{x}\right),E\right),
\]
and the absorption opacity (cross-section) $\sigma_{\varepsilon}^{a}\left(T\left(\mathbf{x}\right),E\right)$
rapidly varies in energy $E$ on the microscale $0<\varepsilon\ll1$;
here $\varepsilon$ denotes the characteristic spacing between spectral
lines or resonances. For simplicity, our discussion ignores density
and pressure dependence in the opacity $\sigma_{\varepsilon}$, but
its incorporation into the proposed algorithm is straightforward;
in fact, see Section~\ref{sub:An-atmospheric-TRT} for an atmospheric
TRT example where the pressure and density dependence are included.
We note that, for neutron problems, $\sigma^{s}$ can vary rapidly
on the micro-scale as well, but we assume her for simplicity that
$\sigma^{s}\left(T,E\right)$ smoothly depends on $E$ (i.e., is independent
of $\varepsilon$). 

The method is based on a rigorous homogenization theory for (\ref{eq:transport equation, intro}),
and has a computational cost that scales independently of the microscale
parameter $\varepsilon$, aside from a pre-computation step analogous
to computing Planck-weighted or Rosseland-weighted opacities. In this
framework, the fine-scale solution $\psi_{\varepsilon}$ converges
to the homogenized solution $\psi_{0}$ in the sense that 
\[
\lim_{\varepsilon_{k}\rightarrow0}\int_{E_{0}}^{E_{1}}\psi_{\varepsilon_{k}}\left(\mathbf{x},\boldsymbol{\Omega},E\right)dE=\int_{E_{0}}^{E_{1}}\psi_{0}\left(\mathbf{x},\boldsymbol{\Omega},E\right)dE,
\]
for arbitrary energy values $E_{0}$ and $E_{1}$ and some sequence
$\varepsilon_{k}\rightarrow0$. Our derivation assumes that the total
opacity $\sigma_{\varepsilon}\left(E,T\right)$ depends on a free
parameter $\varepsilon>0$ (governing the characteristic spacing between
lines), and that $\left\Vert \sigma_{\varepsilon}\left(\cdot,T\right)\right\Vert _{\infty}$
is uniformly bounded in $\varepsilon>0$ for each $T$. In practice,
the final algorithm only makes use of a single opacity $\sigma_{\varepsilon}\left(T,E\right)$
at some fixed characteristic microscale $\varepsilon$ that is much
smaller than the macroscopic scale we are interested in capturing.

As a key tool, we use the Young measure $\lambda_{E}^{T}$ associated
with $\sigma_{\varepsilon}\left(E,T\right)$ (cf. \cite{BALL:1989}).
Roughly speaking, $\lambda_{E}^{T}$ gives the probability distribution
of values $\sigma_{\varepsilon}\left(E,T\right)$ in a vanishingly
small neighborhood of $E$ as $\varepsilon\rightarrow0$. The key
property that the Young measure $\lambda_{E}^{T}$ satisfies (see
\cite{BALL:1989} for a proof) is that, for any continuous function
$F\left(E,\xi\right)$ defined for $E\geq0$ and $\xi\geq0$, there
is a sequence $\varepsilon_{k}\rightarrow0$ such that 
\begin{equation}
\lim_{\varepsilon_{k}\rightarrow0}\int_{E_{0}}^{E_{1}}F\left(E,\sigma_{\varepsilon_{k}}\left(E,T\right)\right)dE=\int_{E_{0}}^{E_{1}}\left(\int_{0}^{\infty}F\left(E,\xi\right)\lambda_{E}^{T}\left(d\xi\right)\right)dE.\label{eq:Young measure, key rep-1}
\end{equation}
Intuitively, for small $\varepsilon$ and for all $E'\in\left[E-\Delta E/2,E+\Delta E/2\right]$
in a neighborhood of $E$ with $0<\varepsilon\ll\Delta E$, the values
$F\left(E',\sigma_{\varepsilon}\left(E'\right)\right)\approx F\left(E,\sigma_{\varepsilon}\left(E'\right)\right)$
can be wildly varying since $\sigma_{\varepsilon}\left(E'\right)$
can rapidly oscillate for $E'\in\left[E-\Delta E/2,E+\Delta E/2\right]$;
however, the average value of $F\left(E,\sigma_{\varepsilon}\left(E'\right)\right)$
for $E'\in\left[E-\Delta E/2,E+\Delta E/2\right]$ is given by weighting
$F\left(E,\xi\right)$ against the probability of $\sigma_{\varepsilon}\left(E\right)\in\left[\xi-d\xi/2,\xi+d\xi/2\right]$,
\begin{equation}
\int_{0}^{\infty}F\left(E,\xi\right)\lambda_{E}^{T}\left(d\xi\right).\label{eq:expected value}
\end{equation}
A simple but instructive example is the Elsasser band model \cite{ELSAS:1938},
\begin{equation}
\sigma_{\varepsilon}\left(E\right)=\frac{\cosh\left(\beta\right)+1}{\cosh\left(\beta\right)-\cos\left(2\pi E/\varepsilon\right)},\label{eq:Elsasser band model, intro}
\end{equation}
$\beta>0$, which models an infinite number of Lorenz lines with equal
spacing of $\varepsilon$ and uniform strength; in Section~\ref{sub:Example-1:-the},
we compute $\lambda_{E}\left(\xi\right)$ for (\ref{eq:Elsasser band model, intro})
analytically.  The Young measure for $\sigma_{\varepsilon}\left(E,T\right)$--and
the associated homogenized solution---exists if $\left\Vert \sigma_{\varepsilon}\left(\cdot,T\right)\right\Vert _{\infty}$
is uniformly bounded in $\varepsilon>0$ for each $T$ \cite{BALL:1989}.
Our derivation will be a direct application of (\ref{eq:Young measure, key rep-1}). 

In order to derive the homogenized equations from (\ref{eq:Young measure, key rep-1}),
we assume that the absorption opacity $\sigma_{\varepsilon}^{a}\left(E,T\right)$
is of the form 
\begin{equation}
\sigma_{\varepsilon}^{a}\left(E,T\right)=\chi\left(\sigma_{\varepsilon}\left(E\right),E,T\right),\label{eq:basic opacity assumption}
\end{equation}
where $\chi\left(\sigma,E,T\right)$ is continuous in its first argument
and $\sigma_{\varepsilon}\left(E\right)$ is an appropriate function
(that may or may not be directly related to the original opacity $\sigma_{\varepsilon}^{a}\left(E,T\right)$).
This assumption, in particular, is a generalization of the commonly
assumed assumption in atmospheric TRT calculations that $\sigma_{\varepsilon}^{a}\left(E,T\right)=\chi\left(\sigma_{\varepsilon}\left(E,T_{0}\right),T\right)$,
where $T_{0}$ is a fixed reference temperature (see \cite{QIA-LIO:1992}
for its use in the correlated k-Distribution method). In Section~\ref{sub:Relaxation-of-the},
we relax the assumption in (\ref{eq:basic opacity assumption}) of
an exact equality, and discuss how to numerically compute an approximation
$\sigma_{\varepsilon}^{a}\left(E,T\right)\approx\chi\left(\sigma_{\varepsilon}\left(E,T_{0}\right),E,T\right)$,
for a given reference temperature $T_{0}$, that is optimal in a certain
sense.

Assume that $\sigma_{\varepsilon}^{a}\left(E,T\right)=\chi\left(\sigma_{\varepsilon}\left(E\right),E,T\right)$.
Then from the solution of the transport equation for $\Psi\left(\mathbf{x},\boldsymbol{\Omega},E,\xi\right)$,
parameterized by the real number $\xi\geq0$,
\begin{eqnarray}
\boldsymbol{\Omega}\cdot\nabla_{\mathbf{x}}\Psi+\left[\chi\left(\xi,E,T\left(\mathbf{x}\right)\right)+\sigma^{s}\left(E,T\left(\mathbf{x}\right)\right)\right]\Psi & = & \chi\left(\xi,E,T\left(\mathbf{x}\right)\right)S\left(T\left(\mathbf{x}\right),E\right)+Q_{0}\left(\mathbf{x},\boldsymbol{\Omega},E\right)+\nonumber \\
 &  & \int_{\mathbb{S}^{2}}\int_{0}^{\infty}\int_{0}^{\infty}\Sigma^{s}\left(\mathbf{x},E,E',\boldsymbol{\Omega}\cdot\boldsymbol{\Omega}'\right)\Psi\,\lambda_{E'}\left(d\xi\right)dE'd\boldsymbol{\Omega}',\label{eq:transport equation on extended phase space, intro}
\end{eqnarray}
we can exactly recover the homogenized solution by weighting against
the Young measure of $\sigma_{\varepsilon}\left(E\right)$, 
\begin{equation}
\psi_{0}\left(\mathbf{x},\boldsymbol{\Omega},E\right)=\int_{0}^{\infty}\Psi\left(\mathbf{x},\boldsymbol{\Omega},E,\xi\right)\lambda_{E}\left(d\xi\right).\label{eq:psi_0 as weighting against Young measure, intro}
\end{equation}
Heuristically, $\Psi\left(\mathbf{x},\boldsymbol{\Omega},E,\xi\right)d\lambda_{E}\left(\xi\right)$
weights the solution of the transport equation (\ref{eq:transport equation on extended phase space, intro})
by the probability that $\sigma_{\varepsilon}\left(E\right)=\xi$
when $\varepsilon\ll1$, and integration over all possible values
$\sigma$ yields the homogenized solution at $\left(\mathbf{x},\boldsymbol{\Omega},E\right)$
in the limit of $\varepsilon\rightarrow0$. The proof of this is a
direct application of (\ref{eq:Young measure, key rep-1}), and is
given in Section~\ref{sec:The-homogenized-transport}. In the proof,
we assume that (\ref{eq:transport equation on extended phase space, intro})
has a unique solution that is continuous in its last argument $\xi$. 

By discretizing equations (\ref{eq:transport equation on extended phase space, intro})
and (\ref{eq:psi_0 as weighting against Young measure, intro}), we
obtain an algorithm that is analogous to the multiband method (cf.
\cite{CUL-POM:1980}). In particular, we choose a coarse number $m_{g}$
of energy groups $\left[E_{i},E_{i+1}\right]$ and a coarse number
$m_{\sigma}$ opacity bands $\left[\sigma_{j},\sigma_{j+1}\right]$,
and use the theory in \cite{BALL:1989} to construct a discrete approximation
\begin{equation}
\lambda_{E}\left(\xi\right)\approx\sum_{j=1}^{m}p_{i,j}\delta\left(\xi-\sigma_{i,j}\right),\,\,\,\, E\in\left[E_{i},E_{i+1}\right],\label{eq:discrete approx of Young measure, intro}
\end{equation}
where $p_{i,j}$ gives the probability that $\sigma_{j}\leq\sigma_{\varepsilon}\left(E\right)\leq\sigma_{j+1}$
for $E\in\left[E_{i},E_{i+1}\right]$ (see Section~\ref{sub:Discrete-approximation-of}
for more details). Convergence of (\ref{eq:discrete approx of Young measure, intro})
to $\lambda_{E}$ is made precise in \cite{BALL:1989}. The key point
of this construction is that, for realistic opacities, the number
parameters $\sigma_{i,j}$ and $p_{i,j}$ is typically a small constant
independent of the size $\varepsilon$ of the microscale; this is
provably so when the opacity is multiscale, e.g. of the form $\sigma_{0}\left(T,E,E/\varepsilon\right)$,
where, e.g., $\sigma_{0}\left(T,E,\kappa\right)$ is almost periodic
in $\kappa$. In general, computing (\ref{eq:discrete approx of Young measure, intro})
scales linearly in $\varepsilon^{-1}$, but needs to be performed
only once for a fine enough energy grid in order for interpolation
to be accurate; this is analogous to the pre-computation of Planck-weighted
or Rosseland-weighted opacities for use in multigroup transport codes,
which also scales linearly in $\varepsilon^{-1}$. 

Now define $\kappa_{ij}\left(T\right)=\chi\left(\sigma_{i,j},E_{i},T\right)$
and 
\[
\Psi_{ij}\left(\mathbf{x},\boldsymbol{\Omega}\right)=\int_{E_{i}}^{E_{i+1}}\Psi\left(\mathbf{x},\boldsymbol{\Omega},E,\sigma_{i,j}\right)dE,\,\,\,\,\,\,\,\, S_{i}\left(E,T\right)=\int_{E_{i}}^{E_{i+1}}S\left(E,T\right)dE.
\]
Then using the representation (\ref{eq:discrete approx of Young measure, intro})
in (\ref{eq:transport equation on extended phase space, intro}),
we obtain the multigroup-type equations 
\begin{eqnarray}
\boldsymbol{\Omega}\cdot\nabla_{\mathbf{x}}\Psi_{ij}+\kappa_{ij}\left(T\left(\mathbf{x}\right)\right)\Psi_{ij} & = & \kappa_{ij}\left(T\left(\mathbf{x}\right)\right)S_{i}\left(x\right)+\nonumber \\
 &  & \int_{\mathbb{S}^{2}}\Sigma^{s}\left(\mathbf{x},E\rightarrow E',\boldsymbol{\Omega}\cdot\boldsymbol{\Omega}'\right)\Psi_{ij}d\boldsymbol{\Omega}',\label{eq:transport multiband-multigroup, intro}
\end{eqnarray}
This corresponds to a standard multigroup-type approximation of (\ref{eq:transport equation on extended phase space, intro})
at each band parameter $\xi=\sigma_{i,j}$. The homogenized solution
at $I_{0}\left(\mathbf{x},\boldsymbol{\Omega},E\right)$ in each group
interval $\left[E_{i},E_{i+1}\right]$ is then given by weighting
against the discrete approximation (\ref{eq:discrete approx of Young measure, intro})
of the Young measure $\lambda_{E}$, 
\begin{equation}
\psi_{0}\left(\mathbf{x},\boldsymbol{\Omega},E\right)\approx\sum_{j=1}^{m}p_{i,j}\Psi_{ij}\left(\mathbf{x},\boldsymbol{\Omega}\right),\,\,\,\, E\in\left[E_{i},E_{i+1}\right].\label{eq:formula for psi_o, discretized Young measures}
\end{equation}
Here we used the discrete approximation (\ref{eq:discrete approx of Young measure, intro})
of the Young measure in (\ref{eq:psi_0 as weighting against Young measure, intro}).
We remark that the values of $p_{i,j}$ can be obtained via interpolation
from a pre-computed table, as explained in Section~\ref{sub:Discrete-approximation-of}.

In our numerical experiments, we choose the band parameters $\sigma_{i,j}$
to be either equally spaced or logarithmically spaced between the
minimum and maximum opacity values in each group, 
\[
\min_{E_{i}\leq E\leq E_{i+1}}\sigma_{\varepsilon}\left(E\right),\,\,\,\,\max_{E_{i}\leq E\leq E_{i+1}}\sigma_{\varepsilon}\left(E\right).
\]
This choice for the band parameters $\sigma_{i,j}$ is justified in
Section~\ref{sub:Discrete-approximation-of}.

To summarize: we solve the multigroup-type equations (\ref{eq:transport multiband-multigroup, intro})
for each group interval $\left[E_{i},E_{i+1}\right]$ and for each
opacity band $\left[\sigma_{i,j},\sigma_{i,j+1}\right]$, and then
average with respect to the discrete Young measure, (\ref{eq:formula for psi_o, discretized Young measures}).
The opacity bands $\left[\sigma_{i,j},\sigma_{i,j+1}\right]$ can
be chosen to be equally spaced or log-spaced within the range of $\sigma_{\varepsilon}\left(E\right)$,
$E\in\left[E_{i},E_{i+1}\right]$.

Equations (\ref{eq:transport multiband-multigroup, intro}) and (\ref{eq:formula for psi_o, discretized Young measures})
are analogous to the multiband method, but are derived from within
a homogenization framework. Unlike the multiband method, however,
no closure assumption (i.e., weighting spectrum) is needed to compute
the multiband parameters $\sigma_{i,j}$. In addition, the group average
(\ref{eq:formula for psi_o, discretized Young measures}) is not a
direct sum as in the multiband method, but instead uses the discrete
Young measure to weight the multiband solutions within each group.
As previously remarked, the homogenization approach is also related
to the correlated k-Distribution method \cite{ARK-GRO:1972} and to
the multigroup k-Distribution method \cite{MOD-ZHA:2001}. We point
out that the current approach is able to explicitly handle scattering
in energy, which to the best of our knowledge has not been explored
with k-Distribution methods.

The remainder of this paper is organized as follows. In Section~\ref{sec:The-homogenized-transport},
we derive the homogenized equations, and their discrete approximation.
We then apply we apply the above methodology to three examples in
Section~\ref{sec:Examples}. For simplicity, we neglect scattering
in all examples. In Section~\ref{sub:Example-1:-the}, we consider
the Elsasser band model (\ref{eq:Elsasser band model, intro}) with
$\varepsilon=10^{-4}$ (see also \cite{POMRAN:1971}). For the Elsasser
band model, we also analytically compute the Young measure $\lambda_{E}\left(\xi\right)$
and compare it with its discrete approximation (\ref{eq:discrete approx of Young measure, intro}).
In the second example, we consider a neutron transport example using
an absoption cross section $\sigma_{\varepsilon}\left(E\right)$ for
iron at room temperature and a Watt fission spectrum; note that, in
this case, the subscript $\varepsilon$ of $\sigma_{\varepsilon}\left(E\right)$
is formerly retained in order to denote the characteristic resonance
spacing, but is not an actual free parameter. Finally, in our last
example, we consider an atmospheric TRT calculation; we use $12$
homogenous atmoshperic layers between $1\,$km-$12\,$km (the temperature
and pressure in each layer come from the 1976 US standard atmosphere),
and consider the absorption opacity corresponding to water vapor.
For all three examples, we demonstrate a small number of energy discretization
parameters can capture the solution with between $0.1$ and $1$ percent
accuracy, and using orders of magnitude fewer parameters than the
standard multigroup formulation with Planck-weighted opacities for
comparable accuracies. We note that, for these examples, we expect
that the correlated k-Distribution method can yield the same accuracy
with a comparable number of parameters.

\section{The homogenized transport equation and its discrete approximation\label{sec:The-homogenized-transport}}

In this Section~\ref{sub:Derivation-of-homogenized}, we first derive
the homogenized equations (\ref{eq:transport equation on extended phase space, intro})
and (\ref{eq:psi_0 as weighting against Young measure, intro}). We
then discuss the discrete approximation of the Young measure in Section~\ref{sub:Discrete-approximation-of},
and derive the discrete approximation (\ref{eq:transport multiband-multigroup, intro})
and (\ref{eq:formula for psi_o, discretized Young measures}) in Section~\ref{sub:Derivation-of-discrete}.
Finally, we conclude this section with a heuristic derivation of the
homogenized system for the time-dependent TRT equations.

\subsection{Derivation of the homogenized system \label{sub:Derivation-of-homogenized}}

To derive (\ref{eq:transport equation on extended phase space, intro})
and (\ref{eq:psi_0 as weighting against Young measure, intro}), first
assume that the scattering kernel only depends on angle, i.e. $\Sigma^{s}=\Sigma^{s}\left(\mathbf{x},\boldsymbol{\Omega}\cdot\boldsymbol{\Omega}'\right)$.
Consider the solution $\Psi\left(\mathbf{x},\boldsymbol{\Omega},E,\xi\right)$
of (\ref{eq:transport equation on extended phase space, intro}).
Then $\Psi\left(\mathbf{x},\boldsymbol{\Omega},E,\sigma_{\varepsilon}\left(E\right)\right)$
satisfies the transport equation (\ref{eq:transport equation, intro}).
Therefore, since (\ref{eq:transport equation, intro}) has a unique
solution, $\psi_{\varepsilon}\left(\mathbf{x},\boldsymbol{\Omega},E\right)=\Psi\left(\mathbf{x},\boldsymbol{\Omega},E,\sigma_{\varepsilon}\left(E\right)\right)$.
It follows from (\ref{eq:Young measure, key rep-1}) that, for any
$0<E_{0}<E_{1}$, 
\begin{eqnarray*}
\lim_{\varepsilon\rightarrow0}\int_{E_{0}}^{E_{1}}\psi_{\varepsilon}\left(\mathbf{x},\boldsymbol{\Omega},E\right)dE & = & \lim_{\varepsilon\rightarrow0}\int_{E_{0}}^{E_{1}}\Psi\left(\mathbf{x},\boldsymbol{\Omega},E,\sigma_{\varepsilon}\left(E\right)\right)dE.\\
 & = & \int_{E_{0}}^{E_{1}}\int_{0}^{\infty}\Psi\left(\mathbf{x},\boldsymbol{\Omega},E,\xi\right)d\lambda_{E}\left(\xi\right)dE
\end{eqnarray*}
Since $E_{0}$ and $E_{1}$ are arbitrary, we finally conclude that
\[
\psi_{0}\left(\mathbf{x},\boldsymbol{\Omega},E\right)=\int_{0}^{\infty}\Psi\left(\mathbf{x},\boldsymbol{\Omega},E,\xi\right)\lambda_{E}\left(d\xi\right).
\]

Now suppose that the scattering kernel $\Sigma^{s}=\Sigma^{s}\left(\mathbf{x},E,E',\boldsymbol{\Omega}\cdot\boldsymbol{\Omega}'\right)$
depends on both energy and angle. We again argue that (\ref{eq:transport equation on extended phase space, intro})
and (\ref{eq:psi_0 as weighting against Young measure, intro}) are
the appropriate homogenized equations. To do so, define $\tilde{\psi}_{\varepsilon}\left(\mathbf{x},\boldsymbol{\Omega},E\right)\equiv\Psi\left(\mathbf{x},\boldsymbol{\Omega},E,\sigma_{\varepsilon}\left(E\right)\right)$.
Then 
\[
\boldsymbol{\Omega}\cdot\nabla_{\mathbf{x}}\tilde{\psi}_{\varepsilon}+\sigma_{\varepsilon}\tilde{\psi}_{\varepsilon}=\sigma_{\varepsilon}^{a}S+Q_{0}+\int_{\mathbb{S}^{2}}\left(\int_{0}^{\infty}\Sigma^{s}\left(\mathbf{x},E,E',\boldsymbol{\Omega}\cdot\boldsymbol{\Omega}'\right)\tilde{\psi}_{\varepsilon}dE'd\boldsymbol{\Omega}'+\mathcal{E}_{\varepsilon}\right)d\boldsymbol{\Omega}',
\]
where the residual term $\mathcal{E}_{\varepsilon}$ is given by 
\[
\mathcal{E}_{\varepsilon}\left(\mathbf{x},\boldsymbol{\Omega}',E\right)=\int_{0}^{\infty}\Sigma^{s}\left(\mathbf{x},E,E',\boldsymbol{\Omega}\cdot\boldsymbol{\Omega}'\right)\left(\tilde{\psi}_{\varepsilon}\left(\mathbf{x},\boldsymbol{\Omega}',E'\right)-\int_{0}^{\infty}\Psi\left(\mathbf{x},\boldsymbol{\Omega}',E',\xi\right)\lambda_{E'}\left(d\xi\right)\right)dE'.
\]
From the property (\ref{eq:Young measure, key rep-1}) and $\tilde{\psi}_{\varepsilon}\left(\mathbf{x},\boldsymbol{\Omega},E'\right)=\Psi\left(\mathbf{x},\boldsymbol{\Omega},E',\sigma_{\varepsilon}\left(E'\right)\right)$,
we see that $\mathcal{E}_{\varepsilon}\left(\mathbf{x},\boldsymbol{\Omega},E\right)\rightarrow0$
for $\varepsilon\rightarrow0$. In addition, the energy dependence
in $\mathcal{E}_{\varepsilon}\left(\mathbf{x},\boldsymbol{\Omega},E\right)$
is only through $\Sigma^{s}\left(\mathbf{x},E,E',\boldsymbol{\Omega}\cdot\boldsymbol{\Omega}'\right)$,
and is therefore slow (i.e., it is independent of the small parameter
$\varepsilon$). It follows that for $\varepsilon\ll1$, 
\[
\boldsymbol{\Omega}\cdot\nabla_{\mathbf{x}}\tilde{\psi}_{\varepsilon}+\sigma_{\varepsilon}\tilde{\psi}_{\varepsilon}-\sigma_{\varepsilon}^{a}S+Q_{0}+\int_{0}^{\infty}\Sigma^{s}\left(\mathbf{x},E,E',\boldsymbol{\Omega}\cdot\boldsymbol{\Omega}'\right)\tilde{\psi}_{\varepsilon}dE\approx0.
\]
Since the transport equation is well-posed, $\tilde{\psi}_{\varepsilon}\approx\psi_{\varepsilon}$.
Finally, invoking (\ref{eq:Young measure, key rep-1}) again, 
\[
\lim_{\varepsilon\rightarrow0}\int_{E_{0}}^{E_{1}}\tilde{\psi}_{\varepsilon}\left(\mathbf{x},\boldsymbol{\Omega},E\right)dE=\int_{E_{0}}^{E_{1}}\int_{0}^{\infty}\Psi\left(\mathbf{x},\boldsymbol{\Omega},E,\xi\right)\lambda_{E}\left(d\xi\right)dE,
\]
and we see that (\ref{eq:transport equation on extended phase space, intro})
and (\ref{eq:psi_0 as weighting against Young measure, intro}) are
the appropriate homogenized equations.

\subsection{Discrete approximation of the Young measures\label{sub:Discrete-approximation-of}}

Here we discuss the derivation of the discrete approximation (\ref{eq:discrete approx of Young measure, intro}),
from which equations (\ref{eq:transport multiband-multigroup, intro})
and (\ref{eq:formula for psi_o, discretized Young measures}) follow
from the homogenized equations (\ref{eq:transport equation on extended phase space, intro})
and (\ref{eq:psi_0 as weighting against Young measure, intro}). As
we will see, the number of bands $\left[\sigma_{j},\sigma_{j+1}\right]$
needed in the approximation to determined by the particular function
$F\left(\xi\right)$ used in the fundamental representation (\ref{eq:Young measure, key rep-1})
for the Young measure; in our applications, $F\left(\xi\right)$ is
very smooth, and a small number of bands are required, independent
of the scale $\varepsilon$ at which $\sigma_{\varepsilon}\left(E\right)$
varies in energy.

We will construct a discrete approximation to the Young measure via
the theory developed in \cite{BALL:1989}; for notational simplicity,
in this section we drop the temperature dependence in the notation.
In particular, the measure $\lambda_{E}$ is entirely determined by
its action on continuous functions $f$ via
\[
\left\langle \lambda_{E},f\right\rangle =\lim_{\Delta E\rightarrow0}\lim_{\varepsilon\rightarrow0}\frac{1}{\Delta E}\int_{E-\Delta E/2}^{E+\Delta E/2}f\left(E',\sigma_{\varepsilon}\left(E'\right)\right)dE'.
\]
Note the order of the limits intuitively corresponds to choosing a
scale $\delta$ that is large relative to the microscopic behavior
but small relative to the macroscopic behavior, $0<\varepsilon\ll\Delta E$.
Then 
\begin{equation}
\left\langle \lambda_{E},f\right\rangle \approx\frac{1}{\Delta E}\int_{E-\Delta E/2}^{E+\Delta E/2}f\left(E,\sigma_{\varepsilon}\left(E'\right)\right)dE',\label{eq:Ball, approximation of Young measure}
\end{equation}
and this becomes precise by first letting $\varepsilon\rightarrow0$
and then $\Delta E\rightarrow0$. In (\ref{eq:Ball, approximation of Young measure}),
approximated $f\left(E',\sigma_{\varepsilon}\left(E'\right)\right)\approx f\left(E,\sigma_{\varepsilon}\left(E'\right)\right)$,
which is valid over $\left[E-\Delta E/2,E+\Delta E/2\right]$ since
the variation of $f$ in its first argument does not depend on the
fast scale $\varepsilon$.

Given fixed $\sigma_{1}<\dots\sigma_{j}<\sigma_{j+1}<\ldots$, define
the characteristic functions $\zeta_{j}\left(\xi\right)$, 
\[
\zeta_{j}\left(\xi\right)=\begin{cases}
1, & \,\sigma_{j}\leq\xi\leq\sigma_{j+1}\\
0, & \,\text{else}.
\end{cases}.
\]
Consider the collection of step functions 
\begin{equation}
f\left(E,\xi\right)=\sum_{j}f\left(E,\sigma_{j}\right)\zeta_{j}\left(\xi\right).\label{eq:step function}
\end{equation}
Although $\zeta_{j}$ are not continuous, any continuous function
can be approximated by such functions. Now, for a general step function
(\ref{eq:step function}) and using (\ref{eq:Ball, approximation of Young measure}),
\begin{eqnarray*}
\left\langle \lambda_{E},f\right\rangle  & \approx & \frac{1}{\Delta E}\int_{E-\Delta E/2}^{E+\Delta E/2}\left(\sum_{j}f\left(E',\sigma_{j}\right)\zeta_{j}\left(\sigma_{\varepsilon}\left(E'\right)\right)\right)dE'\\
 & \approx & \sum_{j}f\left(E,\sigma_{j}\right)\left(\frac{1}{\Delta E}\int_{E-\Delta E/2}^{E+\Delta E/2}\zeta_{j}\left(\sigma_{\varepsilon}\left(E'\right)\right)dE'\right)\\
 & = & \sum_{j}f\left(E,\sigma_{j}\right)p_{j}\left(E\right).
\end{eqnarray*}
In the last equality, the probability $p_{j}\left(E\right)$ is given
by 
\begin{equation}
p_{j}\left(E\right)=\frac{\lambda\left(E'\in\left[E-\Delta E/2,E+\Delta E/2\right]\mid\sigma_{j}\leq\sigma_{\varepsilon}\left(E'\right)\leq\sigma_{j+1}\right)}{\Delta E},\label{eq:definition of p_j(E)}
\end{equation}
with $\lambda$ in (\ref{eq:definition of p_j(E)}) denoting the Lebesque
measure. Approximating a general continuous function $f\left(\xi\right)$
(defined for $\xi\geq0$) by a step function, we have that 
\begin{eqnarray*}
\left\langle \lambda_{E},f\right\rangle  & \approx & \sum_{j}p_{j}\left(E\right)f\left(E,\sigma_{j}\right)\\
 & = & \int f\left(E,\xi\right)\left(\sum_{j}p_{j}\left(E\right)\delta\left(\xi-\sigma_{j}\right)\right)d\xi.
\end{eqnarray*}
Thus, a discrete approximation to the Young measure is given by 
\[
\lambda_{E}\left(\xi\right)\approx\sum_{j}p_{j}\left(E\right)\delta\left(\xi-\sigma_{j}\right),
\]
where $p_{j}\left(E\right)$ is defined via (\ref{eq:definition of p_j(E)}). 

The function $p_{j}\left(E\right)$ give the probability that $\sigma_{j}\leq\sigma_{\varepsilon}\left(E'\right)\leq\sigma_{j+1}$
for $E'\in\left[E-\Delta E/2,E+\Delta E/2\right]$, where $0<\varepsilon\ll\Delta E$;
in particular, $p_{j}\left(E\right)$ may be computed by uniformly
sampling $E'\in\left[E-\Delta E/2,E+\Delta E/2\right]$ and counting
how many times $\sigma_{j}\leq\sigma_{\varepsilon}\left(E'\right)\leq\sigma_{j+1}$
for each band $\left[\sigma_{j},\sigma_{j+1}\right]$. We remark that
$p_{j}\left(E\right)$ can be precomputed on a fine energy grid, and
evaluated at other energy points via interpolation.

\subsection{Derivation of the discrete homogenized system\label{sub:Derivation-of-discrete}}

To derive (\ref{eq:transport multiband-multigroup, intro}), evaluate
(\ref{eq:transport equation on extended phase space, intro}) at $\xi=\sigma_{i,j}$
and integrate in $E$ over $\left[E_{i},E_{i+1}\right]$. Now, approximate
the integral
\[
\int_{E_{i}}^{E_{i+1}}\chi\left(\sigma_{i,j},E,T\left(\mathbf{x}\right)\right)\Psi\left(\mathbf{x},\boldsymbol{\Omega},E,\sigma_{i,j}\right)dE
\]
by
\[
\chi\left(\sigma_{i,j},E_{i},T\left(\mathbf{x}\right)\right)\int_{E_{i}}^{E_{i+1}}\Psi\left(\mathbf{x},\boldsymbol{\Omega},E,\sigma_{i,j}\right)\Psi dE=\kappa_{ij}\left(T\left(\mathbf{x}\right)\right)\Psi_{ij}\left(\mathbf{x},\boldsymbol{\Omega},E\right).
\]
This step is accurate since, by assumption, $\chi\left(\xi,E,T\left(\mathbf{x}\right)\right)$
smoothly varies in energy $E$. We perform a similar calculation for
the integral of $\chi\left(\sigma_{i,j},E,T\left(\mathbf{x}\right)\right)S\left(T\left(\mathbf{x}\right),E\right)$
over $\left[E_{i},E_{i+1}\right]$, and obtain the multigroup-type
system (\ref{eq:transport multiband-multigroup, intro}).

To derive (\ref{eq:formula for psi_o, discretized Young measures}),
we use the discrete approximation (\ref{eq:discrete approx of Young measure, intro})
in (\ref{eq:psi_0 as weighting against Young measure, intro}),
\begin{eqnarray*}
\psi_{0}\left(\mathbf{x},\boldsymbol{\Omega},E\right) & = & \int_{0}^{\infty}\Psi\left(\mathbf{x},\boldsymbol{\Omega},E,\xi\right)\lambda_{E}\left(d\xi\right)\\
 & \approx & \int_{0}^{\infty}\Psi\left(\mathbf{x},\boldsymbol{\Omega},E,\xi\right)\left(\sum_{j=1}^{m}p_{i,j}\delta\left(\xi-\sigma_{i,j}\right)\right)d\xi\\
 & = & \sum_{j=1}^{m}p_{i,j}\Psi\left(\mathbf{x},\boldsymbol{\Omega},E,\sigma_{i,j}\right).
\end{eqnarray*}
Finally, integrating both sides of
\[
\psi_{0}\left(\mathbf{x},\boldsymbol{\Omega},E\right)=\sum_{j=1}^{m}p_{i,j}\Psi\left(\mathbf{x},\boldsymbol{\Omega},E,\sigma_{i,j}\right)
\]
over $\left[E_{i},E_{i+1}\right]$, we obtain (\ref{eq:formula for psi_o, discretized Young measures}).

\subsection{Relaxation of the correlated opacity assumption\label{sub:Relaxation-of-the}}

Our derivation of the homogenized system assumed that the correlated
assumption (\ref{eq:basic opacity assumption}) is an equality for
appropriate functions $\chi\left(x,E,T\right)$ and $\sigma_{\varepsilon}\left(E\right)$.
In this section, we assume instead that 
\begin{equation}
\sigma_{\varepsilon}\left(E,T\right)\approx\chi\left(\sigma_{\varepsilon}\left(E,T_{0}\right),E,T\right),\label{eq:correlated assumption, relaxation section}
\end{equation}
\emph{approximately} holds, where $T_{0}$ is some appropriate reference
temperature. That is, we approximate the opacity at a general temperature
$T$ as functionally related to the opacity as some reference temperature
$T_{0}$. As discussed below, we also include the possibility of slow
energy variation in the functional relationship (\ref{eq:correlated assumption, relaxation section}),
which arises naturally in our following discussion. In general, the
correlated assumption (\ref{eq:correlated assumption, relaxation section})
does not hold. We therefore discuss in this section how to compute
a function $\chi\left(x,E,T\right)$ that best approximates $\sigma_{\varepsilon}\left(E,T\right)\approx\chi\left(\sigma_{\varepsilon}\left(E,T_{0}\right),E,T\right)$.

To motivate the basic idea, consider two random variables $X\geq0$
and $Y\geq0$ and the associated joint probability density $p\left(x,y\right)=\mathbb{P}\left(X=x,Y=y\right)$.
Then it is a well-known fact that the conditional expected value,
\begin{equation}
\chi\left(x\right)\equiv\int_{0}^{\infty}yp\left(x,y\right)dy,\label{eq:conditional expected value}
\end{equation}
minimizes the mean squared error 
\[
\mathbb{E}\left(Y-\chi\left(X\right)\right)^{2}=\min_{g}\mathbb{E}\left(Y-g\left(X\right)\right)^{2}\equiv\int_{0}^{\infty}\int_{0}^{\infty}\left(y-g\left(x\right)\right)^{2}p\left(x,y\right)dxdy,
\]
among all ($X$-measurable) functions $g$. In other words, $\chi\left(X\right)$
is the best functional fit to $Y$ in the sense of minimizing the
average mean squared error.

This naturally leads us to consider the joint Young measure $\lambda_{E,T}\left(\xi_{1},\xi_{2}\right)$
associated with the pair of functions $\left(\sigma_{\varepsilon}\left(E,T_{0}\right),\sigma_{\varepsilon}\left(E,T\right)\right)$
(we explicitly include the temperature $T$ in the notation $\lambda_{E,T}$
to emphasize this dependence). Heuristically, $\lambda_{E,T}\left(\xi_{1},\xi_{2}\right)$
gives the probability density that $\sigma_{\varepsilon}\left(E',T_{0}\right)=\xi_{1}$
and $\sigma_{\varepsilon}\left(E',T\right)=\xi_{2}$ for $E'$ in
a small neighborhood of $E-\delta\leq E'\leq E+\delta$, where $\delta$
is large relative to the characteristic line spacing $\varepsilon$
but small relative to the macroscopic variation of interest (i.e.,
$\varepsilon\ll\delta\ll1$). For small $\varepsilon$, we then have
from (\ref{eq:conditional expected value}) that the conditional expected
value, 
\begin{equation}
\chi\left(x,E,T\right)\equiv\int_{0}^{\infty}\xi_{2}\lambda_{E,T}\left(x,d\xi_{2}\right),\label{eq: chi, conditional expected value}
\end{equation}
approximately minimizes the average error, 
\[
\mathbb{E}\left(\sigma_{\varepsilon}\left(E,T\right)-\chi\left(\sigma_{\varepsilon}\left(E,T_{0}\right),E,T\right)\right)^{2}\equiv\int_{0}^{\infty}\int_{0}^{\infty}\left(\xi_{2}-\chi\left(\xi_{1}\right)\right)^{2}\lambda_{E,T}\left(d\xi_{1},d\xi_{2}\right),
\]
in the limit of small $\varepsilon$. Note that, if the assumption
(\ref{eq:correlated assumption, relaxation section}) exactly holds,
then the joint probability measure $\lambda_{E}\left(\xi_{1},\xi_{2}\right)$
is supported on the curve $\xi_{2}=\chi\left(\xi_{1}\right)$. 

To numerically approximate (\ref{eq: chi, conditional expected value}),
suppose that $\Delta E$ is chosen so that $0<\varepsilon\ll\Delta E$.
Divide the range of $\sigma_{\varepsilon}\left(E',T_{0}\right)$ and
$\sigma_{\varepsilon}\left(E',T\right)$, $E'\in\left[E-\Delta E/2,E+\Delta E/2\right]$,
into temperature-dependant and energy-dependant ``bands'' $\sigma_{j}\left(E,T_{0}\right)$
and $\sigma_{j'}\left(E,T\right)$. For example, in Section~\ref{sub:An-atmospheric-TRT},
we take logarithmically spaced bands between the minimum and maximum
opacity values in $\left[E-\Delta E/2,E+\Delta E/2\right]$. Now define
the discrete probabilities 
\begin{equation}
p_{j,k}\left(E,T\right)\approx\frac{\lambda\left(E'\in\left[E-\Delta E/2,E+\delta/2\right]\mid\sigma_{j}\left(E,T\right)\leq\sigma_{\varepsilon}\left(E',T\right)\leq\sigma_{j+1}\left(E,T\right),\sigma_{j'}\left(E,T_{0}\right)\leq\sigma_{\varepsilon}\left(E',T_{0}\right)\leq\sigma_{j'+1}\left(E,T_{0}\right)\right)}{\delta},\label{eq:p_=00007Bi,j=00007D(E,T)}
\end{equation}
where $\lambda$ again denotes the Lebesgue measure. Then using the
same reasoning as in Section~\ref{sub:Discrete-approximation-of},
we approximate $\chi\left(x,E,T\right)$ as 
\[
\chi\left(x,E,T\right)\approx\sum_{j'}p_{j,j'}\left(E,T\right)\sigma_{j'}\left(E,T\right),\,\,\,\,\text{if}\,\,\,\,\sigma_{j}\left(E,T_{0}\right)\leq x<\sigma_{j+1}\left(E,T_{0}\right).
\]
In practice, we compute the discrete probabilities $p_{j,j'}\left(E,T\right)$
by uniformly sampling energy values $E'$ in $\left[E-\Delta E/2,E+\Delta E/2\right]$
and counting the number of samples for which $\sigma_{j}\left(E,T\right)\leq\sigma_{\varepsilon}\left(E',T\right)\leq\sigma_{j+1}\left(E,T\right)$
and $\sigma_{j'}\left(E,T_{0}\right)\leq\sigma_{\varepsilon}\left(E',T_{0}\right)\leq\sigma_{j'+1}\left(E,T_{0}\right)$.

We approximate the joint Young measure by 
\[
\lambda_{E,T}\left(\xi_{1},\xi_{2}\right)\approx\sum_{j,j'}p_{j,j'}\left(E,T\right)\delta\left(\xi_{1}-\sigma_{j}\left(E,T_{0}\right)\right)\delta\left(\xi_{2}-\sigma_{j'}\left(E,T\right)\right).
\]
Then from (\ref{eq: chi, conditional expected value}), 
\begin{eqnarray*}
\chi\left(x,E,T\right) & \approx & \int_{0}^{\infty}\xi_{2}\left(\sum_{j,j'}p_{j,j'}\left(E,T\right)\delta\left(x-\sigma_{j}\left(E,T_{0}\right)\right)\delta\left(\xi_{2}-\sigma_{j'}\left(E,T\right)\right)\right)d\xi_{2}\\
 & = & \sum_{j,j'}\sigma_{j'}\left(E,T\right)p_{j,j'}\left(E,T\right)\delta\left(x-\sigma_{j}\left(E,T_{0}\right)\right).
\end{eqnarray*}
Integrating $x$ from $\sigma_{j}\left(E,T_{0}\right)$ to $\sigma_{j+1}\left(E,T_{0}\right)$,
\[
\int_{\sigma_{j}\left(E,T_{0}\right)}^{\sigma_{j+1}\left(E,T_{0}\right)}\chi\left(x,E,T\right)dx\approx\sum_{j'}\sigma_{j'}\left(E,T\right)p_{j,j'}\left(E,T\right).
\]

To summarize: in the discrete version of the homogenized system (\ref{eq:transport multiband-multigroup, intro}),
we take 
\begin{equation}
\kappa_{ij}\left(T\right)=\chi\left(\sigma_{i,j},E_{i},T\right)\approx\sum_{j'}p_{j,j'}\left(E_{i},T\right)\sigma_{j'}\left(E_{i},T\right),\label{eq:kappa_=00007Bi,j=00007D}
\end{equation}
where $E_{i}$ denotes the left end point of the $i$th coarse group,
the temperature-dependent probabilities $p_{j,k}\left(E_{i},T\right)$
are defined by (\ref{eq:p_=00007Bi,j=00007D(E,T)}), and the temperature-dependent
bands $\sigma_{k}\left(E_{i},T\right)$ are, e.g., logarithmically
spaced between $\min_{E\in\left[E_{i},E_{i+1}\right]}\sigma_{\varepsilon}\left(E,T\right)$
and $\max_{E\in\left[E_{i},E_{i+1}\right]}\sigma_{\varepsilon}\left(E,T\right)$.
Notice that the probabilities $p_{j,j'}\left(E,T\right)$ may be pre-computed
on a fine energy and temperature grid and evaluated at arbitrary energy
and temperature values via interpolation.

\section{Numerical Examples\label{sec:Examples}}

We apply this methodology to three examples. We first consider in
Section~\ref{sub:Example-1:-the} the radiative transfer equation
at constant temperature and using the Elsasser band opacity (\ref{eq:Elsasser band model, intro}),
where we take the line spacing $\varepsilon=10^{-4}$; this example
is also considered in \cite{POMRAN:1971}. For this simple but instructive
example, we can compute the Young measure analytically and compare
it to its discrete approximation. 

In our second example, we consider a neutron transport problem using
the absorption cross section for iron at room temperature and a Watt
fission spectrum for our source. Whereas opacities contain lines,
nuclear cross sections for neutron applications contain resonances,
which are similar to lines. The cross sections in natural iron contain
thousands of fine resonances much like previous examples contained
many lines.

Our final example is an atmospheric TRT calculation using $12$ homogeneous
atmospheric layers from $0-15$ km, and taking a cross-section corresponding
to water vapor over the frequency interval $1000\leq\nu\leq2000$
(in units of 1/cm); the cross-section for water vapor exhibits thousands
of lines in this frequency range.

Let us discuss the approximation scheme first for Sections~\ref{sub:Example-1:-the}~and~\ref{sub:Example-2:-iron},
since the discretization schemes are essentially identical; we discuss
the approximation scheme for the atmospheric problem in more detail
in Section~\ref{sub:An-atmospheric-TRT}. 

We consider a transport equation of the form 
\begin{equation}
\mu\partial_{x}\psi_{\varepsilon}\left(x,\mu,E\right)+\sigma_{\varepsilon}\left(E\right)\psi_{\varepsilon}\left(x,\mu,E\right)=S\left(E\right).\label{eq:transport in slab geometry}
\end{equation}
In Section~\ref{sub:Example-1:-the}, $\sigma_{\varepsilon}\left(E\right)$
denotes the Elasser band opacity (\ref{eq:Elsasser band model, intro})
and $S\left(E\right)$ denotes $\sigma_{\varepsilon}\left(E\right)B\left(E,T\right)$,
with $B\left(E,T\right)$ denoting the Planck function at constant
temperature; in Section~\ref{sub:Example-2:-iron}, $\sigma_{\varepsilon}\left(E\right)$
denotes the cross-section for iron at room temperature and $S\left(E\right)$
denotes a Watts fission spectrum. 

To compute the discrete approximation of the Young measure $\lambda_{E}^{T}$,
we approximate for each energy group $\left[E_{i},E_{i+1}\right]$,
\begin{equation}
\mu_{i}^{T}\left(\sigma\right)\approx\sum_{j=1}^{m}p_{i,j}\delta\left(\sigma-\sigma_{i,j}\right),\label{eq:discrete rep of Young measure, example}
\end{equation}
where $p_{i,j}$ is proportional to the probability that $\sigma_{i,j}\leq\sigma_{\varepsilon}\left(E\right)\leq\sigma_{i,j+1}$
for $E_{i}\leq E\leq E_{i+1}$; $p_{i,j}$ is computed by uniformly
sampling random numbers from $\left[E_{i},E_{i+1}\right]$ (the number
of samples is chosen to be much larger than the number of energy values
needed to resolve $\sigma_{\varepsilon}\left(E\right)$ in $\left[E_{i},E_{i+1}\right]$),
counting how many times $ $$\sigma_{i,j}\leq\sigma_{\varepsilon}\left(E\right)\leq\sigma_{i,j+1}$
for each band $\left[\sigma_{i,j},\sigma_{i,j+1}\right]$, and normalizing
by the total number of samples. We emphasize that, although evaluating
$p_{i,j}$ scales linearly in $\varepsilon^{-1}$, this is a pre-computation
and need only be done once; this pre-computation is analogous to computing
Planck-weighted or Rosseland-weighted opacities.

Using the discrete representation (\ref{eq:discrete rep of Young measure, example}),
we have that for $\mu>0$, 
\begin{eqnarray}
\lim_{\varepsilon\rightarrow0}\int_{E_{i}}^{E_{i+1}}\psi_{\varepsilon}\left(x,\mu,E\right)dE & = & \int_{E_{i}}^{E_{i+1}}S\left(E\right)\left(\int_{0}^{\mathbb{R}}\frac{\left(1-e^{-x\left(\sigma/\mu\right)}\right)}{\sigma}d\lambda_{E}^{T}\left(\sigma\right)\right)dE\nonumber \\
 & \approx & \int_{E_{i}}^{E_{i+1}}S\left(E\right)\left(\sum_{j=1}^{m}p_{i,j}\frac{\left(1-e^{-x\left(\sigma_{i,j}/\mu\right)}\right)}{\sigma_{i,j}}\right)dE\nonumber \\
 & = & S_{i}\sum_{j=1}^{m}p_{i,j}\frac{\left(1-e^{-x\left(\sigma_{i,j}/\mu\right)}\right)}{\sigma_{i,j}},\label{eq:mu > 0, homogenized}
\end{eqnarray}
where 
\[
S_{i}=\int_{E_{i}}^{E_{i+1}}S\left(E\right)dE.
\]
Similarly, for $\mu<0$,
\begin{equation}
\lim_{\varepsilon\rightarrow0}\int_{E_{i}}^{E_{i+1}}\psi_{\varepsilon}\left(s,E\right)dE\approx S_{i}\sum_{j=1}^{m}p_{i,j}\left(1-e^{-\left(x-1\right)\left(\sigma_{i,j}/\mu\right)}\right).\label{eq:mu < 0, homogenized}
\end{equation}

In Sections~\ref{sub:Example-1:-the}~and~\ref{sub:Example-2:-iron},
we compare the discrete approximation, (\ref{eq:mu > 0, homogenized})-(\ref{eq:mu < 0, homogenized}),
to the exact solution integrated over $\left[E_{i},E_{i+1}\right]$,
\[
\int_{E_{i}}^{E_{i+1}}\psi_{\varepsilon}\left(x,\mu,E\right)dE=\int_{E_{i}}^{E_{i+1}}S\left(E\right)\frac{\left(1-e^{-\left(\sigma_{\varepsilon}\left(E\right)/\mu\right)x}\right)}{\sigma_{\varepsilon}\left(E\right)}dE,
\]
for both the Elsasser band (\ref{eq:Elsasser band model, intro})
and for iron opacity as generated via the NJoy \cite{NJOY} program.
More precisely, we compute the exact energy-integrated solution, 
\[
\int_{E_{i}}^{E_{i+1}}\psi_{\varepsilon}\left(x_{i},\mu_{j},E\right)dE,
\]
for $10$ equispaced spatial points $x_{i}\in\left[0,1\right]$ and
for $8$ Gauss-Legendre nodes $\mu_{j}$, as well as the homogenized
energy-integrated solution. We compare the ``exact'' scalar flux
(that is, exact to within angular discretization errors), 
\begin{equation}
\varphi_{\varepsilon}\left(x_{k}\right)=\sum_{j}\sum_{i}\int_{E_{i}}^{E_{i+1}}\psi_{\varepsilon}\left(x_{k},\mu_{j},E\right)dE,\label{eq:discrete scalar flux, exact}
\end{equation}
against its homogenized version 
\begin{equation}
\varphi_{0}\left(x_{k}\right)=\sum_{j}\sum_{i}\int_{E_{i}}^{E_{i+1}}\psi_{0}\left(x_{k},\mu_{j},E\right)dE.\label{eq:discrete scalar flux, approx}
\end{equation}

We compare the results to the standard multigroup method using Planck-weighted
and Rosseland-weighted opacities. In particular, we integrate (\ref{eq:transport equation, intro})
over $\left[E_{i},E_{i+1}\right]$, 
\[
\mu\partial_{x}\int_{E_{g}}^{E_{g+1}}\psi_{\varepsilon}\left(x,\mu,E\right)dE=\int_{E_{g}}^{E_{g+1}}\sigma_{\varepsilon}\left(E,T\right)\left(\frac{S\left(E\right)}{\sigma_{\varepsilon}\left(E\right)}-\psi_{\varepsilon}\left(x,\mu,E\right)\right)dE.
\]
We write 
\[
\int_{E_{g}}^{E_{g+1}}\sigma_{\varepsilon}\left(E\right)\left(\frac{S\left(E\right)}{\sigma_{\varepsilon}\left(E\right)}-\psi_{\varepsilon}\left(x,\mu,E\right)\right)dE\approx\sigma_{g}\left(\int_{E_{g}}^{E_{g+1}}\frac{S\left(E\right)}{\sigma_{\varepsilon}\left(E\right)}dE-\int_{E_{g}}^{E_{g+1}}\psi_{\varepsilon}\left(x,\mu,E\right)dE\right),
\]
where 
\[
\sigma_{g}=\frac{\int_{E_{g}}^{E_{g+1}}S\left(E\right)dE}{\int_{E_{g}}^{E_{g+1}}S\left(E\right)/\sigma_{\varepsilon}\left(E\right)dE}.
\]
Thus, we need to solve 
\[
\mu\partial_{x}\psi_{g}\left(x,\mu\right)+\sigma_{g}\psi_{g}\left(x,\mu\right)=S_{g}\left(T\right),
\]
where 
\[
S_{g}\left(T\right)=\int_{E_{g}}^{E_{g+1}}\frac{S\left(E\right)}{\sigma_{\varepsilon}\left(E\right)}dE,\,\,\,\,\sigma_{g}=\frac{\int_{E_{g}}^{E_{g+1}}S\left(E\right)dE}{\int_{E_{g}}^{E_{g+1}}S\left(E\right)/\sigma_{\varepsilon}\left(E\right)dE}.
\]

\subsection{A regular band model example\label{sub:Example-1:-the}}

In order to assess the accuracy of the discrete approximation of
the Young measure $\lambda_{E}$ associated with (\ref{eq:Elsasser band model, intro}),
we first compute $\lambda_{E}$ analytically (note that, since (\ref{eq:Elsasser band model, intro})
does not depend on temperature $T$, we drop the superscript $T$
on \textbf{$\lambda_{E}^{T}$}). 

Using that $\sigma_{\varepsilon}\left(E\right)=\sigma_{1}\left(E/\varepsilon\right)$,
with $\sigma_{1}\left(E\right)$ a $1$-periodic function, it is a
standard result (see e.g. \cite{BESICO:1954}) that 
\[
\lim_{\varepsilon\rightarrow0}\int_{E_{0}}^{E_{1}}f\left(\sigma_{\varepsilon}\left(E'\right)\right)dE'=\int_{0}^{1}f\left(\sigma_{1}\left(E'\right)\right)dE',
\]
for any values $0\leq E<E_{1}\leq1$. Now, consider the change of
variables, 
\[
\xi=\frac{\cosh\left(\beta\right)+1}{\cosh\left(\beta\right)-\cos\left(2\pi E\right)}.
\]
Then
\[
d\xi=-\frac{2\pi}{\cosh\left(\beta\right)+1}\xi^{2}\sqrt{1-\frac{\left(\xi\cosh\left(\beta\right)-\cosh\left(\beta\right)-1\right)^{2}}{\xi^{2}}}dE.
\]
It follows that, with $c_{\beta}=\left(\cosh\left(\beta\right)+1\right)/\left(\cosh\left(\beta\right)-1\right)$,
\[
\int_{0}^{1}f\left(\sigma_{1}\left(E'\right)\right)dE'=\frac{\cosh\left(\beta\right)+1}{2\pi}\int_{1}^{c_{\beta}}\frac{f\left(\xi\right)}{\xi^{2}\sqrt{1-\xi^{-2}\left(\sigma\cosh\left(\beta\right)-\cosh\left(\beta\right)-1\right)^{2}}}d\xi.
\]
Therefore, the Young measure is given by
\begin{equation}
d\lambda_{E}\left(\xi\right)=\frac{1}{\xi^{2}\sqrt{1-\xi^{-2}\left(\sigma\cosh\left(\beta\right)-\cosh\left(\beta\right)-1\right)^{2}}}\chi_{\left[1,c_{\beta}\right]}\left(\xi\right)d\xi.\label{eq:analytic Young measure}
\end{equation}
As expected, $\lambda_{E}\left(\xi\right)$ is independent of energy
$E$.

\begin{figure}
\protect\caption{Comparison of the Young measure (\ref{eq:analytic Young measure})
and its discrete approximation (\ref{eq:discrete approx of Young measure, Elsasser}).
We use $m=30$ ``band boundaries'' $\sigma_{j}$ in (\ref{eq:discrete approx of Young measure, Elsasser}).
\label{fig:young-measure}}

\includegraphics[scale=0.6]{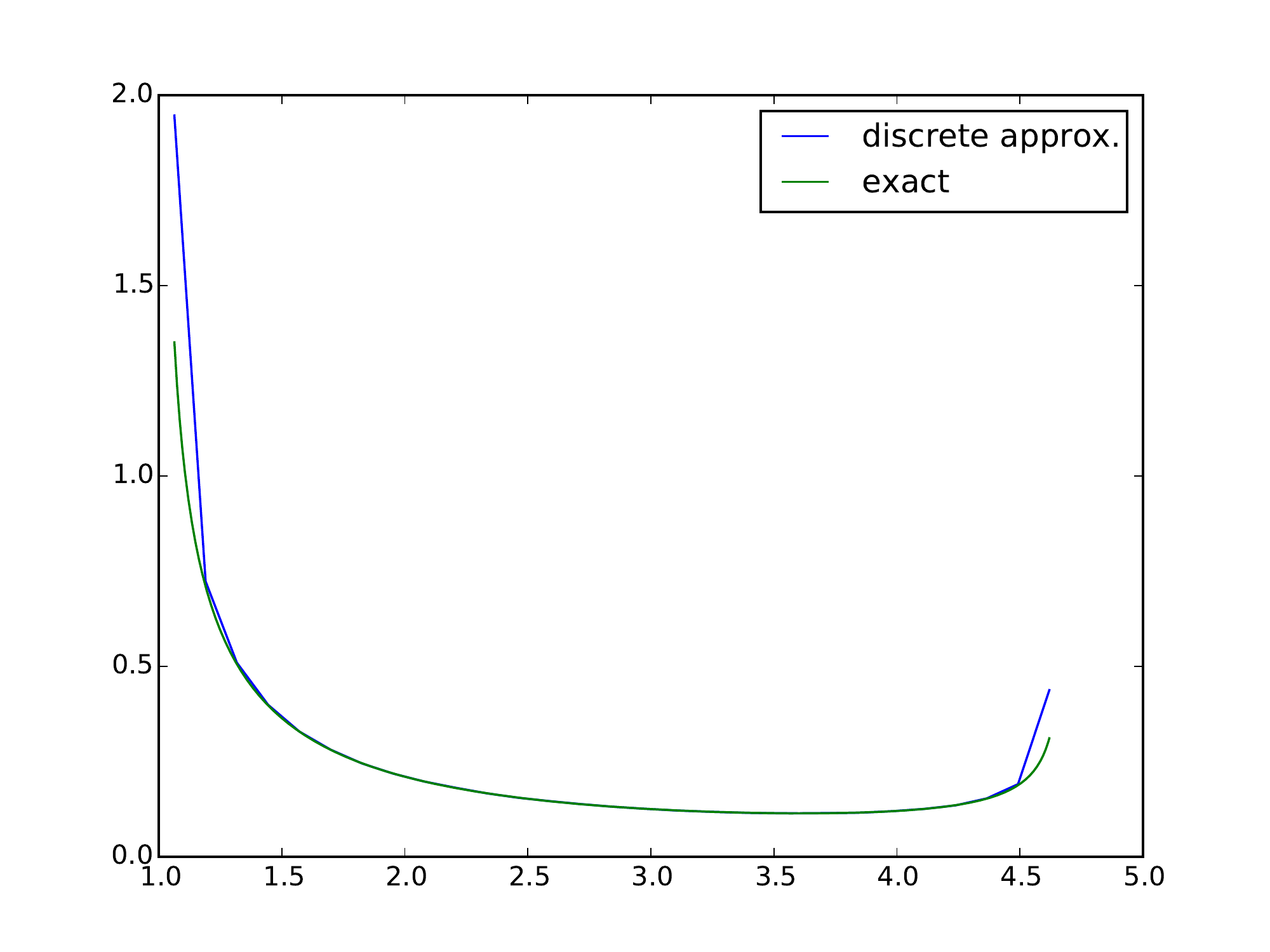}
\end{figure}

We follow the discussion in (\ref{sub:Discrete-approximation-of})
and approximate 
\begin{equation}
\lambda_{E}\left(\xi\right)\approx\sum_{j=1}^{m}p_{j}\delta\left(\xi-\xi_{j}\right).\label{eq:discrete approx of Young measure, Elsasser}
\end{equation}
For this example, $p_{j}$ and $\xi$ are independent of $E$. In
Figure~\ref{fig:young-measure}, we compare the Young measure (\ref{eq:analytic Young measure})
associated with (\ref{eq:Elsasser band model, intro}) against its
discrete approximation (\ref{eq:discrete approx of Young measure, Elsasser});
in this example, we take $\beta=1$ and, in (\ref{eq:discrete approx of Young measure, Elsasser}),
$m=30$ equispaced values $\sigma_{j}$ between $1$ and $c_{\beta}$,
$\beta=1$. Note that the error in the discrete approximation (\ref{eq:discrete approx of Young measure, Elsasser})
necessarily rises at the end points, since the density $d\lambda_{E}\left(\xi\right)$
is infinite at $\sigma=1$ and $\sigma=c_{\beta}$.

Using (\ref{eq:discrete approx of Young measure, Elsasser}) in (\ref{eq:mu > 0, homogenized})-(\ref{eq:mu < 0, homogenized}),
we then obtain the (approximate) homogenized solution $\psi_{0}\left(x,\mu,E\right)$.
In Figure~\ref{fig:transport comparison, Elsasser}, we compare the
``exact'' scalar flux (\ref{eq:discrete scalar flux, exact}) against
its homogenized version (\ref{eq:discrete scalar flux, approx}),
where $\psi_{0}$ is computed from (\ref{eq:discrete approx of Young measure, Elsasser})
and (\ref{eq:mu > 0, homogenized})-(\ref{eq:mu < 0, homogenized}).
We use $m=30$ equispaced values $\sigma_{j}$ between $1$ and $c_{\beta}$,
$\beta=1$, and a single group $\left[E_{0},E_{1}\right]=\left[0,1\right]$. 

~

\begin{figure}
\protect\caption{Comparison of the exact scalar flux (\ref{eq:discrete scalar flux, exact})
against its discrete approximation (\ref{eq:discrete scalar flux, approx}),
using $m=30$ equispaced values $\sigma_{j}$ between $1$ and $c_{\beta}$,
$\beta=1$, and a single group $\left[E_{0},E_{1}\right]=\left[0,1\right]$.
\label{fig:transport comparison, Elsasser}}

\includegraphics[scale=0.6]{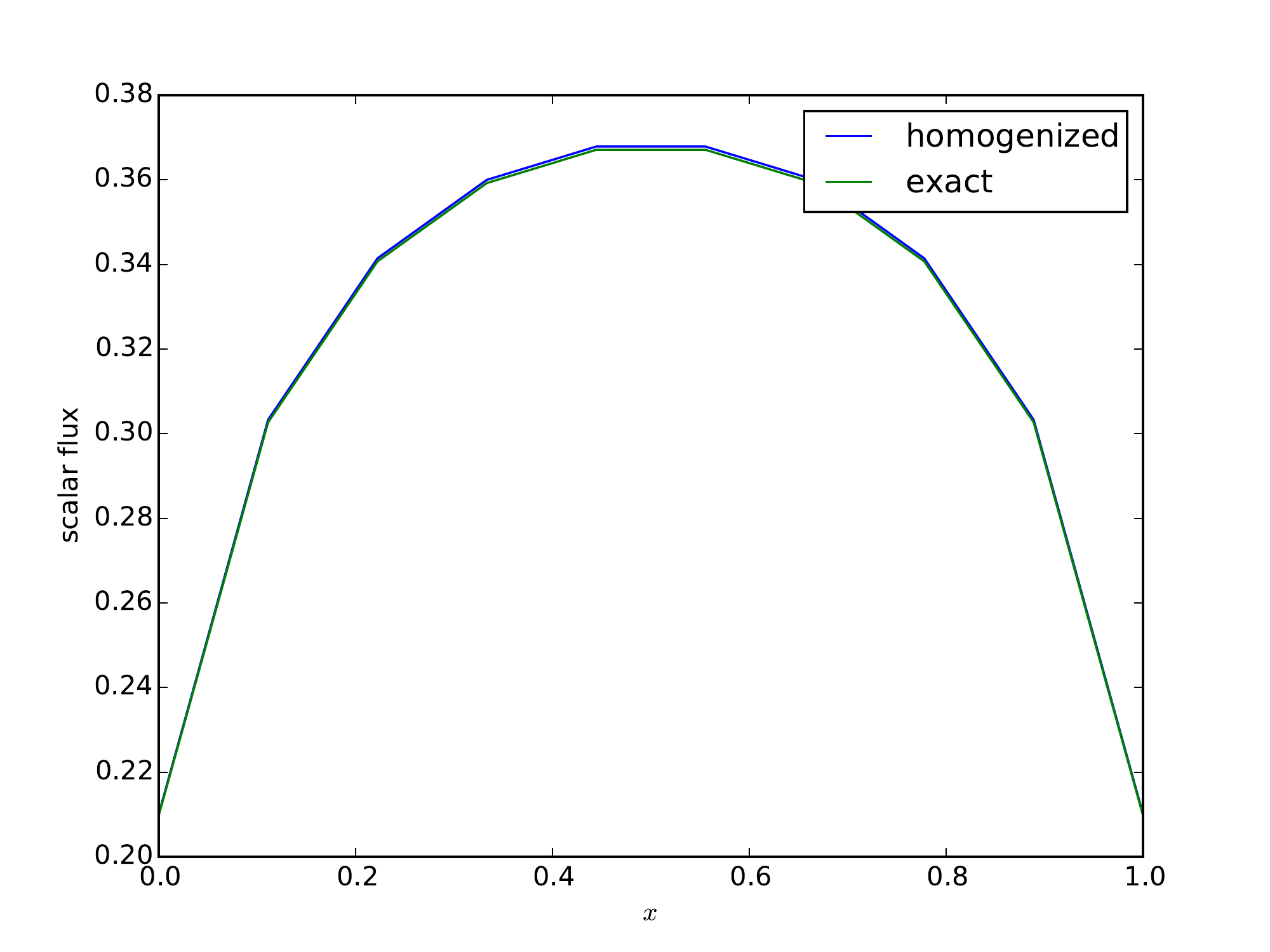}
\end{figure}

\subsection{A neutron transport example with absorption cross section for iron
\label{sub:Example-2:-iron}}

~

\begin{figure}
\protect\caption{Watt fission spectrum}

\includegraphics[scale=0.4]{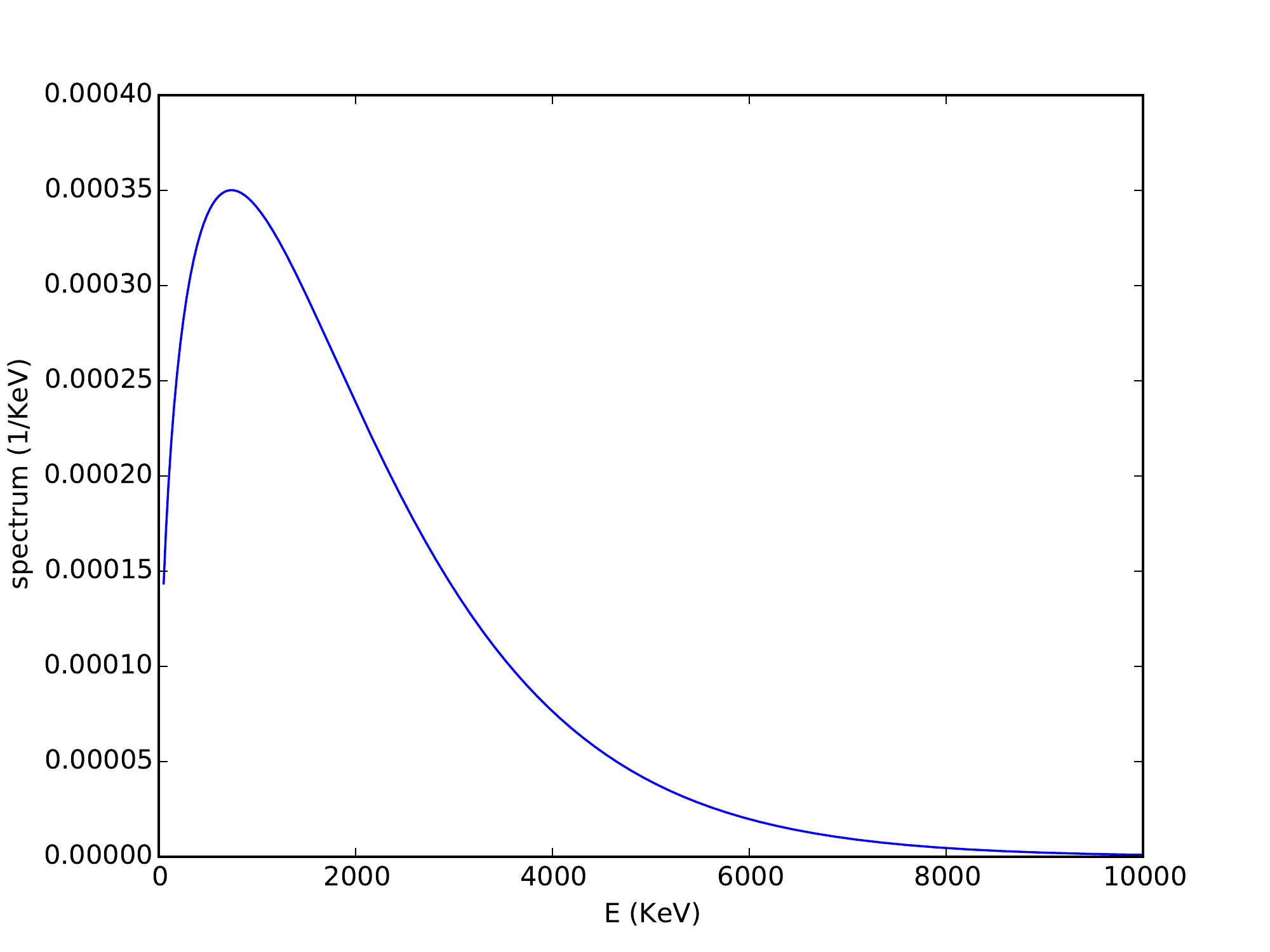}
\end{figure}

~

\begin{figure}
\protect\caption{Absorption cross section for iron.\label{fig:opacity for iron}}

\includegraphics[scale=0.6]{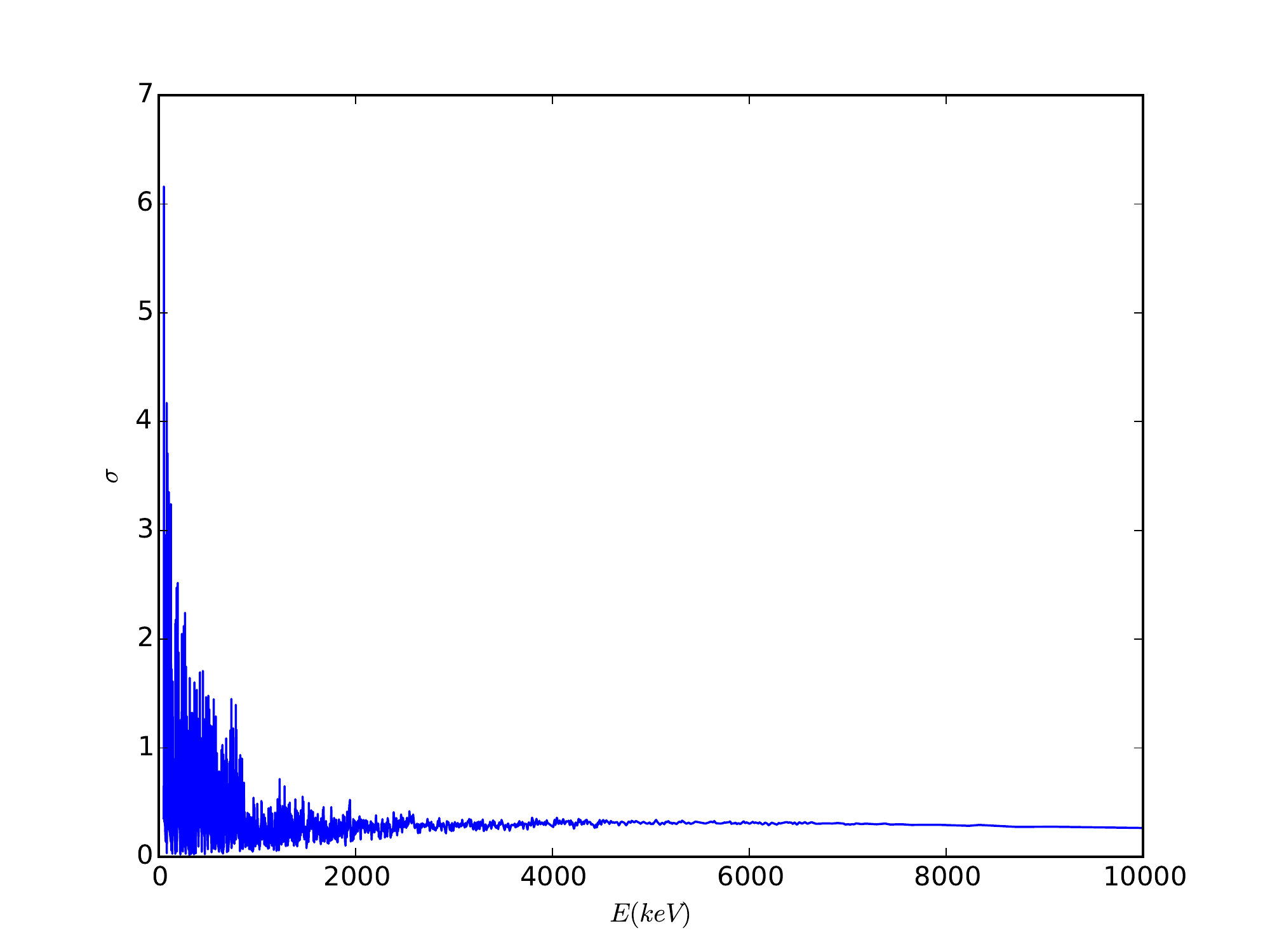}
\end{figure}

\begin{figure}
\protect\caption{Comparison of the exact scalar flux (\ref{eq:discrete scalar flux, exact})
against its homogenized approximation (\ref{eq:discrete scalar flux, approx})
and the multigroup method with Planck-weighted cross sections, using
(a) $1,2024$, (b) $4,096$, (c) $16384$, and (d) $65,536$ equally
spaced groups. Here the cross section is for iron as shown in Figure~\ref{fig:opacity for iron}.
\label{fig:comparison of scalar fluxes, iron}}

(a)

~

\includegraphics[scale=0.4]{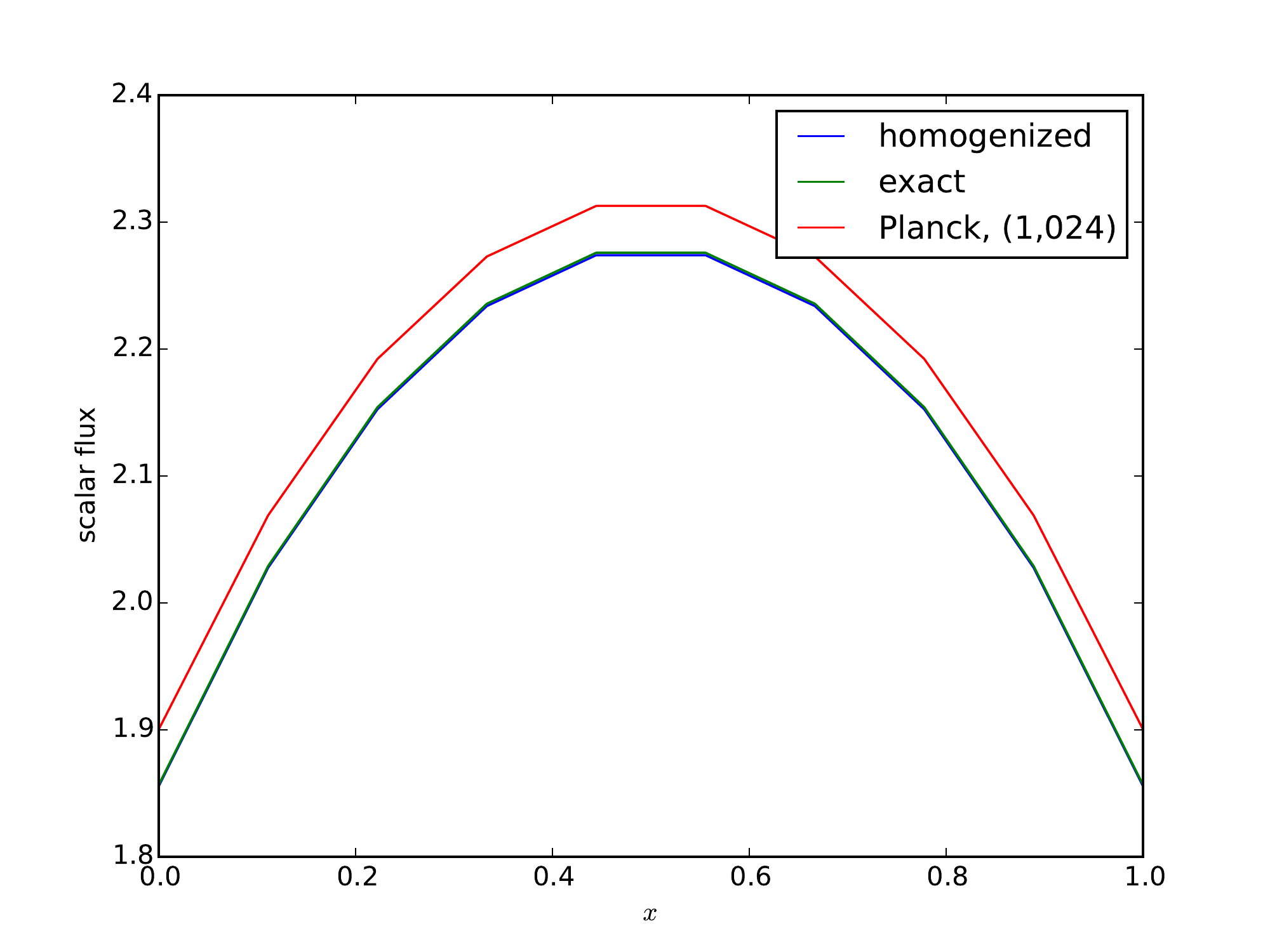}~~~~~~~~

~

(b)

~~

\includegraphics[scale=0.4]{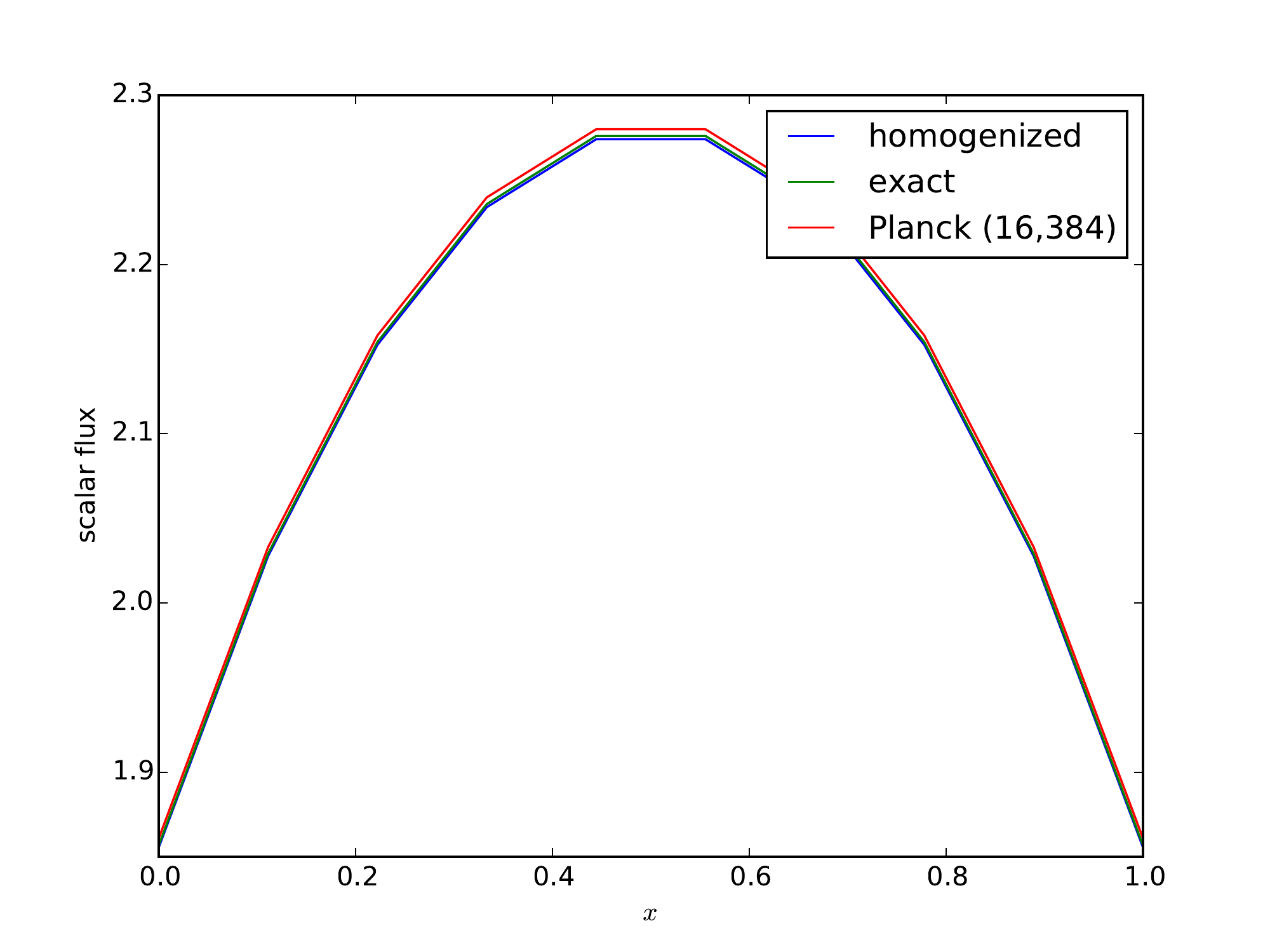}

~

(c)

~~

\includegraphics[scale=0.4]{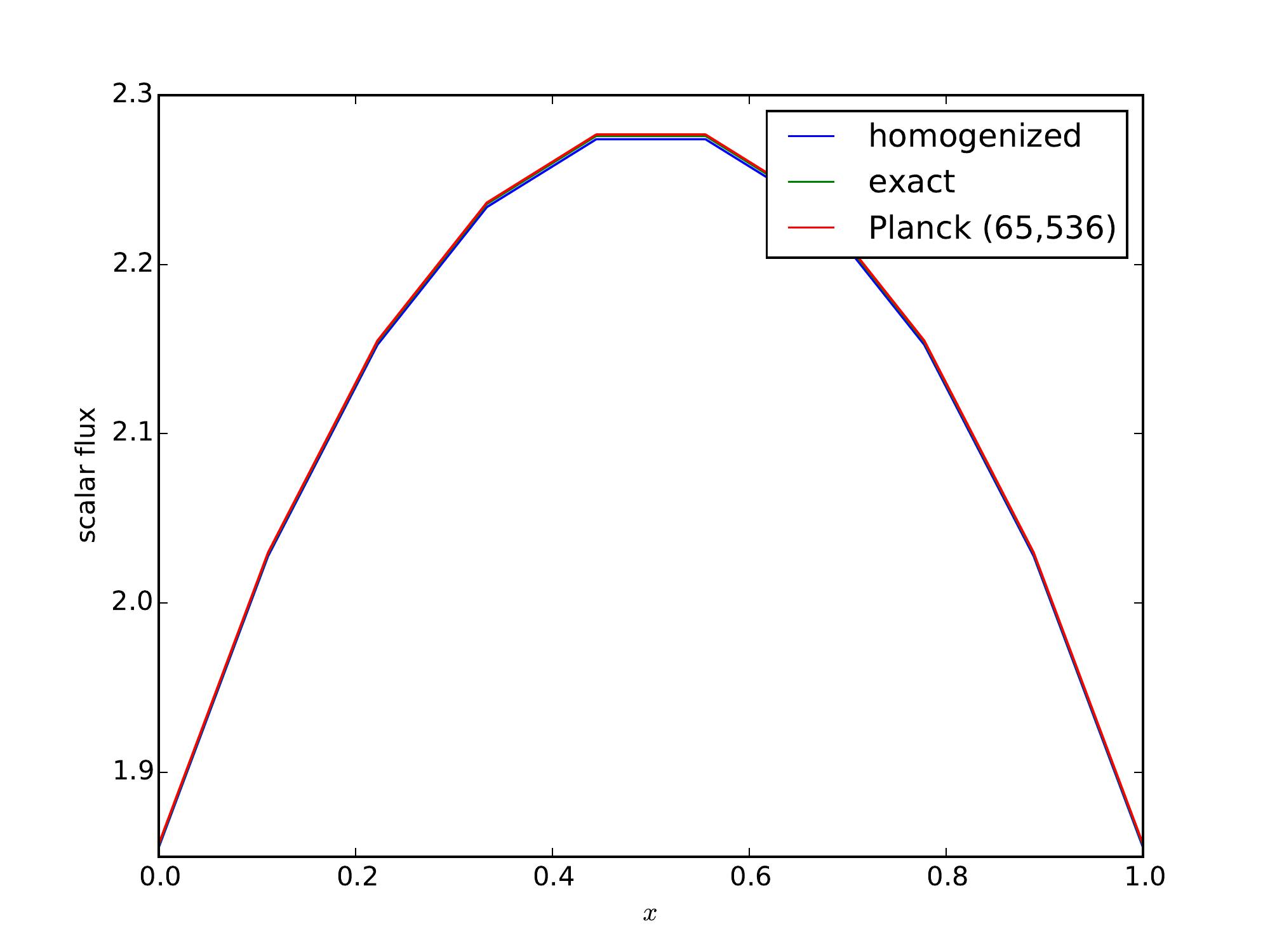}
\end{figure}

\begin{figure}
\protect\caption{Relative error of the homogenized scalar flux (\ref{eq:discrete scalar flux, approx}),
as compared to the exact scalar flux (\ref{eq:discrete scalar flux, exact}),
using the cross section for iron shown in Figure~\ref{fig:opacity for iron}.
\label{fig:error in scalar flux, iron}}

\includegraphics[scale=0.6]{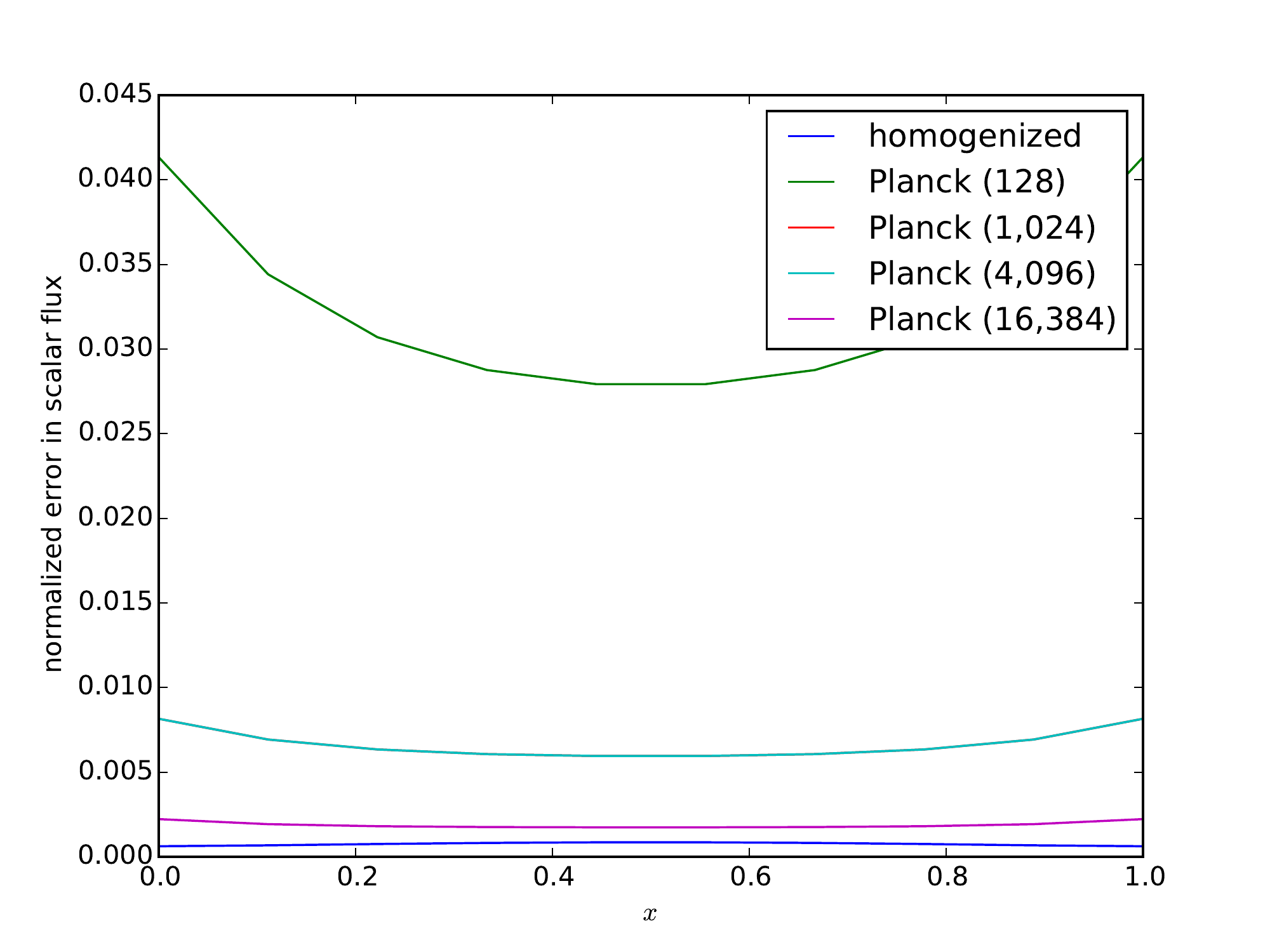}
\end{figure}

We consider a neutron transport example in slab geometry, which describes
the uncollided neutron flux and is useful in shielding applications.
In particular, in equation (\ref{eq:transport in slab geometry}),
we take $\sigma_{\varepsilon}\left(E\right)$ to correspond to the
absorption cross section of iron at room temperature over the energy
range $50$ keV to $10$ MeV, as displayed in Figure~\ref{fig:opacity for iron};
energy/cross-section pairs $\left(E_{j},\sigma_{j}\right)$ for iron
at were generated via the NJoy problem and the cross section $\sigma_{\varepsilon}\left(E\right)$
at a general value $E$ is evaluated via linear interpolation. The
subscript $\varepsilon$ on $\sigma_{\varepsilon}\left(E\right)$
is retained for notational consistency, and represents the characteristic
resonance spacing; note, however, that $\varepsilon$ is not a free
parameter. We use a source term $S\left(E\right)$ corresponding to
a Watt fission spectrum,

\[
\begin{aligned}S(E) & =c\exp\left(-E/a\right)\sinh\sqrt{bE}\end{aligned}
,
\]
with $a=0.988$ MeV, $b=2.2249$ MeV$^{-1}$, and $c=\exp\left(-ab/4\right)/\sqrt{\left(\pi a^{3}b/4\right)}$
MeV$^{-1}$ a normalization constant. 

As in Section~\ref{sub:Example-1:-the}, we compare the exact scalar
flux (\ref{eq:discrete scalar flux, exact}) (that is, exact to within
angular discretization errors) against its homogenized version (\ref{eq:discrete scalar flux, approx}).
We take $4$ energy groups $\left[E_{i},E_{i+1}\right]$, where the
group boundaries $E_{i}$ are equally spaced between $E_{\min}=50.002$
(KeV) and $E_{\max}=10^{4}$ (KeV) Within each energy group $\left[E_{i},E_{i+1}\right]$,
we use $m=40$ equispaced values $\sigma_{j}$ between 
\[
\sigma_{\min,i}=\min_{E_{i}\leq E\leq E_{i+1}}\sigma_{\varepsilon}\left(E\right),
\]
and 
\[
\sigma_{\max,i}=\max_{E_{i}\leq E\leq E_{i+1}}\sigma_{\varepsilon}\left(E\right).
\]
In Figure~\ref{fig:comparison of scalar fluxes, iron}, we compare
the homogenized scalar flux with the exact scalar flux, as well as
with the scalar flux obtained using the multigroup method with $1,024$,
$ $ $16,384$, and $65,536$ equally groups. Figure~\ref{fig:error in scalar flux, iron}
displays the relative errors in the scalar for in the homogenized
scalar flux and the scalar flux using the Planck-weighted opacities;
we see from this figure that the homogenized scalar flux and the multigroup
method $16,384$, achieves about a $.1$ percent relative error in
comparison to the exact scalar flux. Recall that the homogenization
approach used $4\times40=160$ energy discretization parameters, and
so for this accuracy this translates into orders of magnitude fewer
parameters. We note that, for a larger error, the efficiency gain
of the homogenized approach in this example is much less significant.

\subsection{An atmospheric TRT example with water vapor\label{sub:An-atmospheric-TRT}}

\begin{figure}
\protect\caption{Cross-section $\kappa\left(\nu,T,p\right)$ (kg/m\textasciicircum{}2)
for water vapor as a function of frequency (1/cm), with $T=288.15$
(K) and $p=8.9874\times10^{4}$ (Pa).\label{fig:Cross-section--(kg/m^2)}}

\includegraphics[scale=0.6]{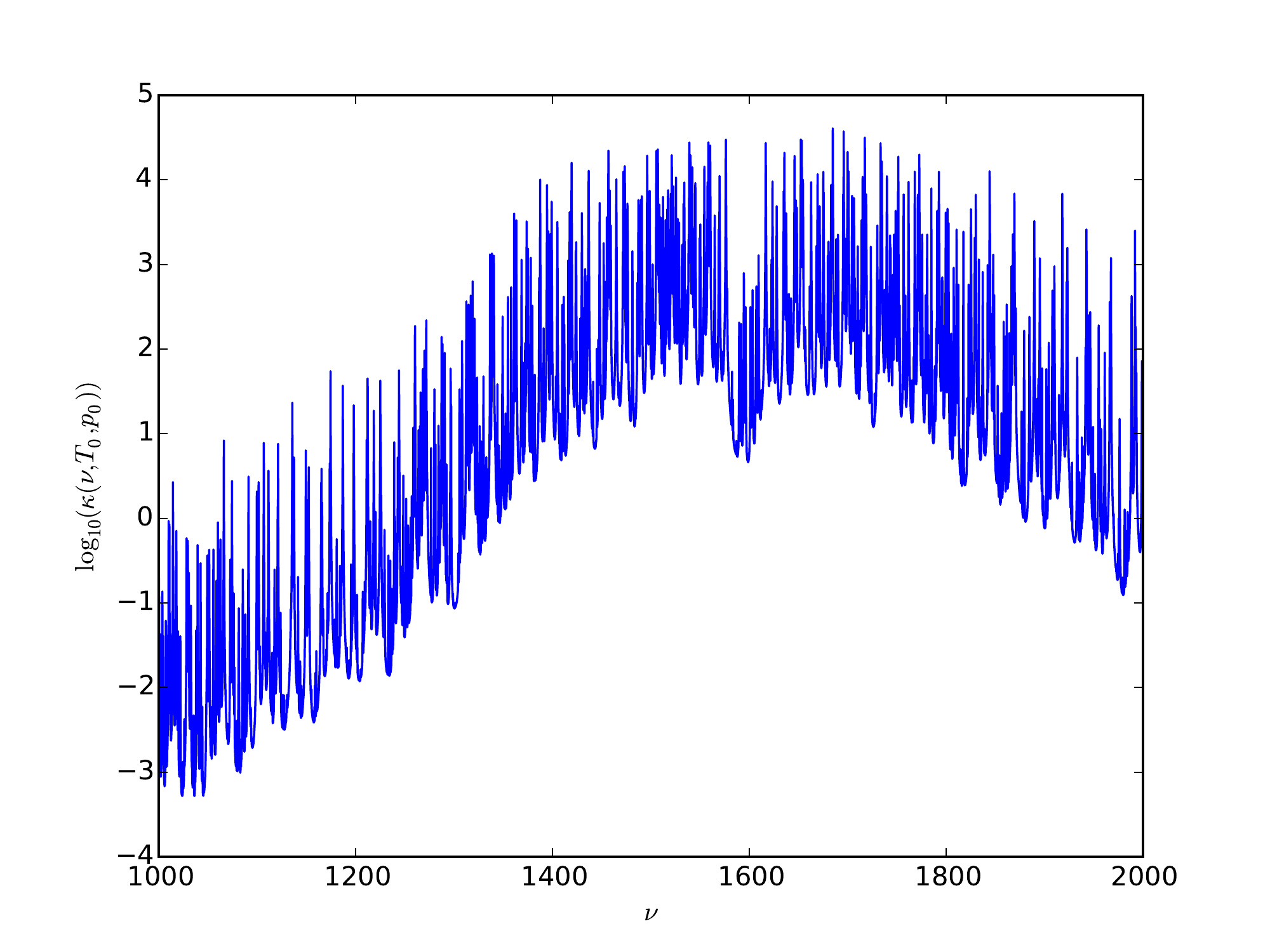}
\end{figure}

\begin{table}
\protect\caption{Temperature and pressure values in each homogeneous layer.\label{tab:Temperature-and-pressure}}

\begin{tabular}{|c|c|c|c|}
\hline 
Height (km) & Temperature (K) & Pressure (Pa) & volume fraction (unitless)\tabularnewline
\hline 
\hline 
0-1 & 281.65 & 8.98746E+4 & .0081\tabularnewline
\hline 
1-2 & 275.15 & 7.94952E+4 & .0077\tabularnewline
\hline 
2-3 & 268.65 & 7.01085E+4 & .0059\tabularnewline
\hline 
3-4 & 262.15 & 6.16402E+4 & .0028\tabularnewline
\hline 
4-5 & 255.65 & 5.40199E+4 & .0016\tabularnewline
\hline 
5-6 & 249.15 & 4.71810E+4 & .0008\tabularnewline
\hline 
6-7 & 242.65 & 4.10607E+4 & .0003\tabularnewline
\hline 
7-8 & 236.15 & 3.55998E+4 & 7.96E-5\tabularnewline
\hline 
8-9 & 229.65 & 3.07425E+4 & 3.21E-5\tabularnewline
\hline 
9-10 & 223.15 & 2.64363E+4 & 1.78E-5\tabularnewline
\hline 
10-12 & 216.65 & 1.93304E+4 & 6.94E-6\tabularnewline
\hline 
12-15 & 216.65 & 1.20446E+4 & 3.84E-6\tabularnewline
\hline 
\end{tabular}
\end{table}

\begin{figure}
\protect\caption{Comparison of the outgoing flux using solutions of the exact and homogenized
equations (\ref{eq:psi, atmosphere}) and (\ref{eq:Psi, atmosphere}).
In solving (\ref{eq:psi, atmosphere}) using the homogenization approach,
we use $10$ energy groups and $5$ bands. We also solve (\ref{eq:psi, atmosphere})
using $2500$ and $5000$ Planck-weighted opacities. Plot (a) shows
the exact $F_{\varepsilon}^{+}\left(x\right)$ and homogenized outgoing
fluxes $F_{0}^{+}\left(x\right)$ (homogenized) and $F_{P}^{+}\left(x\right)$
(Planck). Plot (b) shows the relative errors in the homogenized solution
and the Planck solutions. \label{fig:comparison of line-by-line and homogenized outgoing fluxes}}

(a)

\includegraphics[scale=0.4]{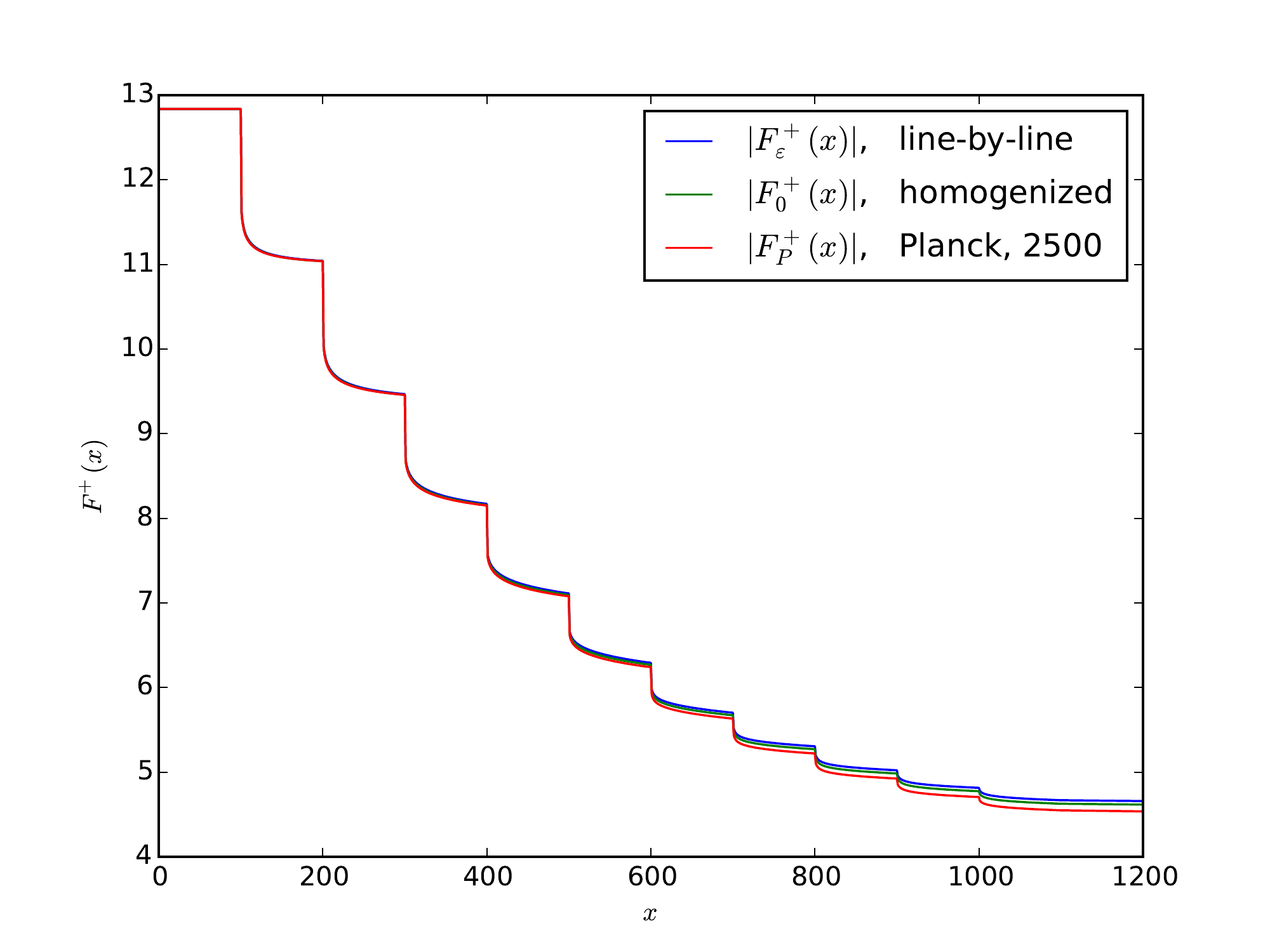}

~

~

(b)

\includegraphics[scale=0.4]{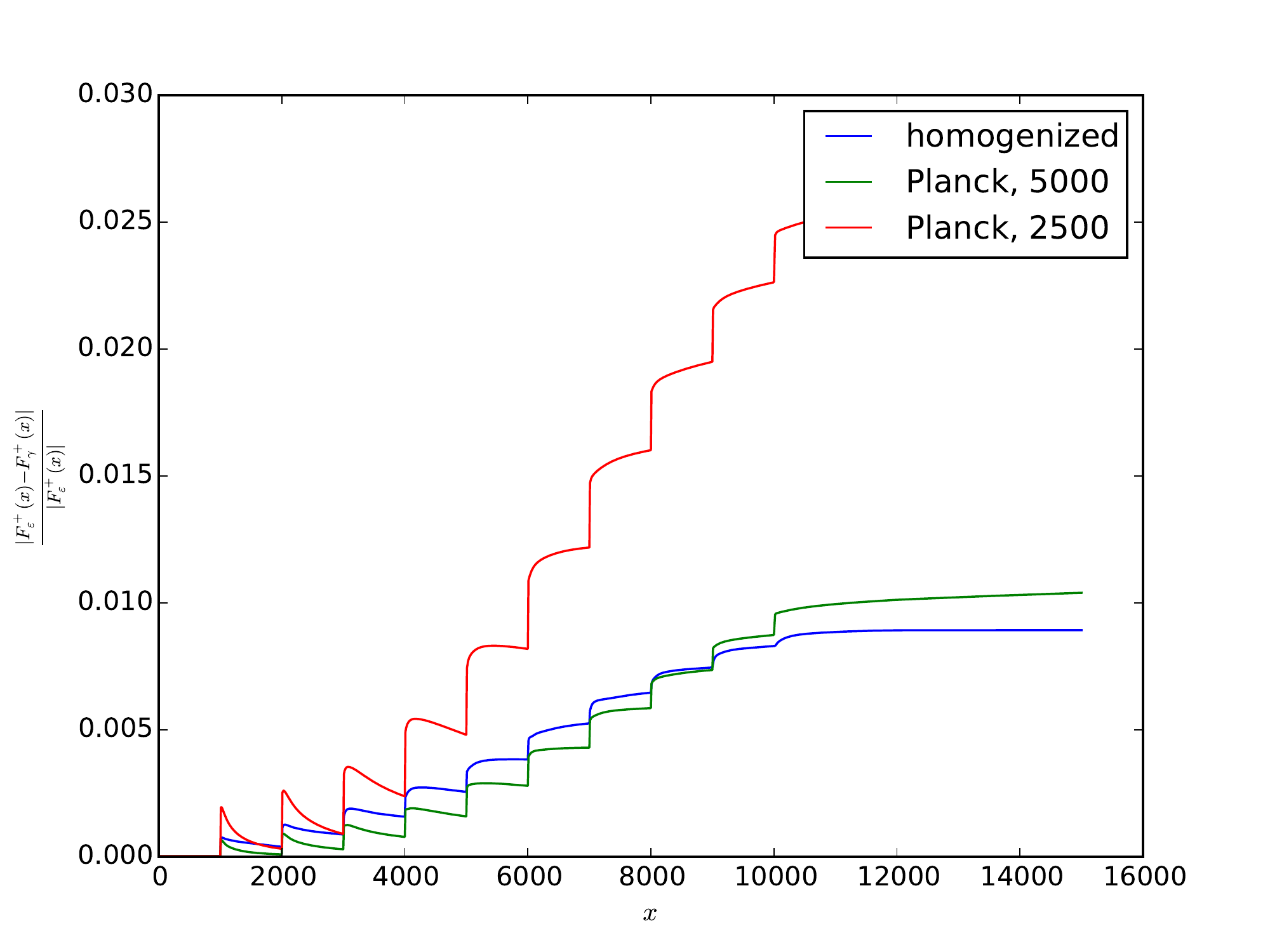}
\end{figure}

We consider an atmospheric model problem consisting of $12$ homogeneous
layers of constant temperature $T\left(x\right)=T_{k}$ and pressure
$p\left(x\right)=p_{k}$ for $x_{k}\leq x\leq x_{k+1}$ (see e.g.
\cite{G-W-C-C:1989}); for consistency with traditional notation,
we consider the frequency variable $\nu$ in place of the energy variable
$E$. We solve the TRT equations using cross-sections $\hat{\sigma}\left(\nu,T,p\right)$,
corresponding to water vapor, over the frequency range $1000-2000$
(1/cm). As standard for infrared atmospheric TRT calculations, we
neglect scattering and assume slab (i.e., plane parallel) geometry.
We also use the 1976 U. S. Standard Atmosphere values for temperature
and pressure in each layer, as displayed in Table~\ref{tab:Temperature-and-pressure}.
In addition, the volume fraction of water vapor in each layer is also
shown in Table~\ref{tab:Temperature-and-pressure}.

The cross-sections $\hat{\sigma}\left(\nu,T,p\right)$ (in units of
kg/m\textasciicircum{}3) are computed using the HITRAN database and
assuming a Lorenz shape profile; in our calculations, we use the Julia
module by J. Bloch-Johnson, https://github.com/jsbj/Jultran.jl/blob/master/src/Jultran.jl,
for processing the HITRAN files in order to generate the opacities.
Here $\log_{10}\left(\hat{\sigma}\left(\nu,T,p\right)\right)$ is
plotted in Figure~\ref{fig:Cross-section--(kg/m^2)} for $T=288.15$
(K) and $p=8.9874\times10^{4}$ (Pa).

We now describe the problem in more detail. In each homogeneous layer
$x\in\left[x_{k},x_{k+1}\right]$, $1\leq k\leq12$, the intensity
$\psi_{\varepsilon}\left(x,\mu,\nu\right)$ satisfies the transport
equation 
\begin{eqnarray}
\mu\partial_{x}\psi_{\varepsilon}^{k}\left(x,\mu,\nu\right)+\sigma_{\varepsilon}\left(\nu,T_{k},p_{k}\right)\psi_{\varepsilon}^{k}\left(x,\mu,\nu\right) & = & \sigma_{\varepsilon}\left(\nu,T_{k},p_{k}\right)B\left(\nu,T_{k}\right),\\
\psi_{\varepsilon}^{k}\left(x_{k},\mu,\nu\right) & = & \psi_{\varepsilon}^{k-1}\left(x_{k},\mu,\nu\right),\,\,\,\,\,\,\,\,\mu>0,\\
\psi_{\varepsilon}^{k}\left(x_{k+1},\mu,\nu\right) & = & \psi_{\varepsilon}^{k+1}\left(x_{k+1},\mu,\nu\right),\,\,\,\,\,\,\,\,\mu<0.\label{eq:psi, atmosphere}
\end{eqnarray}
Here $B\left(\nu,T\right)$ denotes the Planck function; $\sigma_{\varepsilon}\left(\nu,T_{k},p_{k}\right)=r_{k}\rho_{k}\hat{\sigma}_{\varepsilon}\left(\nu,T_{k},p_{k}\right)$,
where $\rho_{k}$ and $r_{k}$ denote the density of air (kg/m\textasciicircum{}3)
and the fraction of water vapor in the $k$th layer and $\hat{\sigma}_{\varepsilon}\left(\nu,T_{k},p_{k}\right)$
denotes the cross-section corresponding to water vapor. We also take
boundary conditions 
\begin{eqnarray*}
\psi_{\varepsilon}\left(x_{0},\mu,\nu\right) & = & B\left(\nu,T_{0}\right),\,\,\,\,\mu>0,\\
\psi_{\varepsilon}\left(x_{K+1},\mu,\nu\right) & = & 0,\,\,\,\,\mu<0.
\end{eqnarray*}

For the homogenized transport equation (\ref{eq:transport multiband-multigroup, intro}),
we use in each spatial interval $\left[x_{k},x_{k+1}\right]$, $k=1,\ldots,12$,
seven logarithmically spaced values $\kappa_{ij}\left(T_{k_{0}}\right)$,
$j=1,\ldots,7$, for each each frequency group $\nu\in\left[\nu_{i},\nu_{i+1}\right]$
and for the reference temperature $T_{k_{0}}$ for layer $k_{0}=6$
(we find that the error is relatively insensitive to this choice);
that is, for each frequency group $\left[\nu_{i},\nu_{i+1}\right]$,
the bands $\log\left(\kappa_{ij}\left(T_{k_{0}}\right)\right)$, $j=1,\ldots,5$,
are equally spaced between the minimum $\min_{\nu_{i}\leq\nu\leq\nu_{i+1}}\log\left(\sigma_{\varepsilon}\left(\nu,T_{k_{0}},p_{k_{0}}\right)\right)$
and the maximum $\max_{\nu_{i}\leq\nu\leq\nu_{i+1}}\log\left(\sigma_{\varepsilon}\left(\nu,T_{k_{0}},p_{k_{0}}\right)\right)$.
For values of $k\neq k_{0}$ (i.e., other spatial intervals), we compute
values $\kappa_{ij}\left(T_{k}\right)$ using equation (\ref{eq:kappa_=00007Bi,j=00007D})
and the techniques discussed in Section~\ref{sub:Relaxation-of-the}. 

From $\kappa_{ij}\left(T_{k}\right)$, we solve 
\begin{eqnarray}
\mu\partial_{x}\Psi_{ij}^{k}\left(x,\mu,\nu\right)+\kappa_{ij}\left(T_{k}\right)\Psi_{ij}^{k}\left(x,\mu,\nu\right) & = & \kappa_{ij}\left(T_{k}\right)B_{i}\left(T_{k}\right),\\
\Psi_{ij}^{k}\left(x_{k},\mu,\nu\right) & = & \Psi_{ij}^{k}\left(x_{k-1},\mu,\nu\right),\,\,\,\,\,\,\,\,\mu>0,\label{eq:Psi, atmosphere}\\
\Psi_{ij}^{k}\left(x_{k},\mu,\nu\right) & = & \Psi_{ij}^{k}\left(x_{k+1},\mu,\nu\right),\,\,\,\,\,\,\,\,\mu<0,
\end{eqnarray}
where 
\[
B_{i}\left(T_{k}\right)=\int_{\nu_{i}}^{\nu_{i+1}}B\left(\nu,T\right).
\]
Given $\Psi_{ij}^{k}\left(x,\mu,\nu\right)$ the homogenized solution
$\psi_{0}$ is computed via 
\[
\psi_{0}\left(x,\mu,\nu\right)=\sum_{j}p_{ij}\Psi_{ij}^{k}\left(x,\mu,\nu\right),\,\,\,\, x\in\left[x_{k},x_{k+1}\right],\,\,\,\nu\in\left[\nu_{i},\nu_{i+1}\right],
\]
where $p_{ij}$ denotes the probability that $\sigma_{j}\leq\sigma_{\varepsilon}\left(\nu,T_{k_{0}},p_{k_{0}}\right)\leq\sigma_{j+1}$
for $\nu_{i}\leq\nu\leq\nu_{i+1}$ (recall that $k_{0}=6$ is taken
for the reference layer). 

In the line-by-line solution of (\ref{eq:psi, atmosphere}), we discretize
in angle using $n_{\mu}=8$ Gaussian quadrature nodes $\mu_{p}$ and
weights $w_{p}$, and $n_{\nu}=200,001$ equally spaced frequency
points in $1000\leq\nu\leq2000$. In each spatial interval $\left[x_{k},x_{k+1}\right]$,
we directly evaluate the analytic solution using $n_{x}=100$ equally
spaced spatial points. For example, for $\mu>0$, we proceed from
the first spatial interval $\left[x_{1},x_{2}\right]$ to the last
spatial interval $\left[x_{11},x_{12}\right]$ and directly evaluate
the analytic solution 
\[
\psi_{\varepsilon}^{k}\left(x,\mu,\nu\right)=e^{-\left(\sigma_{\varepsilon}\left(\nu,T_{k},p_{k}\right)/\mu\right)\left(x-x_{k}\right)}\psi_{\varepsilon}^{k-1}\left(x_{k},\mu,\nu\right)+\left(1-e^{-\left(\sigma_{\varepsilon}\left(\nu,T_{k},p_{k}\right)/\mu\right)\left(x-x_{k}\right)}\right)B\left(\nu,T_{k}\right),
\]
at $n_{x}$ equally spaced points $x_{k}\leq x_{q}^{k}\leq x_{k+1}$,
with $x_{1}^{k}=x_{k}$ and $x_{n_{x}}^{k}=x_{k+1}$. 

Similarly, in the solution of (\ref{eq:Psi, atmosphere}), we discretize
in angle using $n_{\mu}=8$ Gaussian quadrature nodes $\mu_{p}$ and
weights $w_{p}$, and $10$ frequency groups $\left[\nu_{i},\nu_{i+1}\right]$.
In each spatial interval $\left[x_{k},x_{k+1}\right]$ and each energy
group $\left[\left[\nu_{i},\nu_{i+1}\right]\right]$, we use $7$
bands $\kappa_{ij}\left(T_{k}\right)$, that are equally spaced on
a logarithmic scale. As in the solution of (\ref{eq:psi, atmosphere}),
we directly evaluate the analytic solution using $n_{x}=100$ equally
spaced spatial points. For example, for $\mu>0$, we proceed from
the first spatial interval $\left[x_{1},x_{2}\right]$ to the last
spatial interval $\left[x_{12},x_{13}\right]$ and directly evaluate
the analytic solution 
\[
\Psi_{ij}^{k}\left(x,\mu\right)=e^{-\left(\kappa_{ij}\left(T_{k}\right)/\mu\right)\left(x-x_{k}\right)}\Psi_{ij}^{k-1}\left(x_{k},\mu\right)+\left(1-e^{-\left(\kappa_{ij}\left(T_{k}\right)/\mu\right)\left(x-x_{k}\right)}\right)B_{i}\left(T_{k}\right),
\]
at $n_{x}$ equally spaced points $x_{k}\leq x_{q}^{k}\leq x_{k+1}$,
with $x_{1}^{k}=x_{k}$ and $x_{n_{x}}^{k}=x_{k+1}$. 

In Figure~\ref{fig:comparison of line-by-line and homogenized outgoing fluxes},
we compare the line-by-line outgoing flux
\[
F_{\varepsilon}\left(x\right)=\sum_{p=1}^{n_{\mu}}\sum_{i=1}^{n_{\nu}}w_{p}\mu_{p}\left(\nu_{i+1}-\nu_{i}\right)\psi_{\varepsilon}\left(x,\mu_{p},\nu_{i}\right),
\]
against its homogenized version 
\[
F_{0}\left(x\right)=\sum_{p=1}^{n_{\mu}}\sum_{i=1}^{n_{G}}w_{p}\mu_{p}\psi_{0}\left(x,\mu_{p},\nu_{i}\right).
\]
We also compare the exact outgoing flux that against that obtained
using $5000$ Planck-weighted opacities. Plot (a) shows the line-by-line
solution, the solution using $2500$ and $5000$ Planck-weighted frequency
groups, and the homogenized solution using $10$ frequency groups
and $7$ bands . We see from plot (b) in Figure~\ref{fig:comparison of line-by-line and homogenized outgoing fluxes}
that, for a comparable error, the homogenized solution requires about
$70\times$ fewer parameters than the solution obtained via Planck-weighted
frequency groups. 

~~

\bibliographystyle{plain}
\bibliography{opacity_averaging}

\end{document}